\documentclass[aap,MSNbibl,seceqn,dvips]{arximspdf}
\usepackage{graphicx}
\doi{10.1214/13-AAP919} 
\volume{24}
\issue{1}
\pubyear{2014}
\firstpage{150}
\lastpage{197}

\makeatletter
\newcommand{\eqref}[1]{(\ref{#1})}
\newcommand{\cB}{{\mathcal{B}}}
\newcommand{\cG}{\mathcal{G}}
\newcommand{\cC}{{\mathcal{C}}}
\newcommand{\cF}{{\mathcal{F}}}
\newcommand{\cI}{{\mathcal{I}}}
\newcommand{\cL}{{\mathcal{L}}}
\newcommand{\cN}{{\mathcal{N}}}
\newcommand{\cY}{{\mathcal{Y}}}
\newcommand{\R}{{\mathbb{R}}}
\newcommand{\Z}{{\mathbb{Z}}}
\newcommand{\NN}{{\mathbb{N}}}
\newcommand{\Zde}{{\mathbb{Z}^d_e}}
\newcommand{\hxi}{\widehat\xi}
\newcommand{\AV}{\mathrm{AV}}
\newcommand{\LV}{\mathrm{LV}}
\newcommand{\vm}{\mathrm{VM}}
\newcommand{\gv}{\mathrm{GV}}
\newcommand{\tvm}{\mathrm{TV}}
\newcommand{\tzeta}{\tilde\zeta}

\newcommand{\vep}{\varepsilon}
\newcommand{\To}{\Rightarrow}
\newcommand{\zero}{\mathbf{0}}
\newcommand{\one}{\mathbf{1}}

\newtheorem{thmm}{Theorem}[section]
\newtheorem{theorem}{Theorem}

\newtheorem{prop}[thmm]{Proposition}
\newtheorem{lem}[thmm]{Lemma}
\newtheorem{cor}[thmm]{Corollary}
\newproclaim{rem}{Remark}
\newproclaim{example}{Example}
\newtheorem{conj}{Conjecture}
\newproclaim{definition}{Definition}
\makeatother

\begin{document}
\begin{frontmatter}

\title{A complete convergence theorem for voter model perturbations}
\runtitle{Complete convergence theorem}\vspace*{6pt}

\begin{aug}
\author[A]{\fnms{J. Theodore} \snm{Cox}\corref{}\thanksref{t1}\ead[label=e1]{jtcox@syr.edu}}
\and
\author[B]{\fnms{Edwin A.} \snm{Perkins}\thanksref{t2}\ead[label=e3]{perkins@math.ubc.ca}}
\thankstext{t1}{Supported in part by grants from the NSA and NSF.}
\thankstext{t2}{Supported in part by an NSERC Discovery grant.}\vspace*{3pt}
\runauthor{J.~T. Cox and E. A. Perkins}
\affiliation{Syracuse University and University of British Columbia}\vspace*{3pt}
\address[A]{Department of Mathematics\\
Syracuse University\\
Syracuse, New York 13244\\
USA\\
\printead{e1}}

\address[B]{Department of Mathematics\\
University of British Columbia\\
1984 Mathematics Road\\
Vancouver, British Columbia V6T 1Z2\\
Canada\\
\printead{e3}}
\end{aug}

\received{\smonth{8} \syear{2012}}
\revised{\smonth{12} \syear{2012}}

%
\begin{abstract}
We prove a complete convergence theorem for a class of
symmetric voter model perturbations with annihilating
duals. A special case of interest covered by our results
is the stochastic spatial Lotka--Volterra model introduced
by Neuhauser and Pacala [\textit{Ann. Appl. Probab.} \textbf
{9} (1999) 1226--1259]. We also treat two
additional models, the ``affine'' and ``geometric'' voter
models.
\end{abstract}

%
\begin{keyword}[class=AMS]
\kwd[Primary ]{60K35}
\kwd{82C22}
\kwd[; secondary ]{60F99}
\end{keyword}

\begin{keyword}
\kwd{Complete convergence theorem}
\kwd{Lotka--Volterra}
\kwd{interacting particle system}
\kwd{voter model perturbation}
\kwd{annihilating dual}
\end{keyword}
\vspace*{6pt}
\end{frontmatter}

\section{Introduction}\label{secintro}
In our earlier study of voter model perturbations
\mbox{\cite{CDP11,CMP,CP05,CP07,CP08}} we found conditions for survival,
extinction and coexistence for these interacting particle
systems. Our goal here is to show that under additional
conditions it is possible to determine all stationary
distributions and their domains of attraction. We start by
introducing the primary example of this work, a competition
model from~\cite{NP}.

The state of the system at time $t$ is represented by a
spin-flip process $\xi_t$ taking values in
$\{0,1\}^{\Z^d}$. The dynamics will in part be determined by
a fixed probability kernel $p\dvtx\Z^d\to[0,1]$. We assume
throughout that
%
%
\begin{equation}
\label{passump} %
\begin{tabular}{p{280pt}@{}}
$p(0)=0, p(x)$ is
symmetric, irreducible, and has covariance  matrix $\sigma^2$I for some  $\sigma^2\in(0,
\infty).$
\end{tabular}
\end{equation}
For most of our results we will need to assume that $p(x)$ has
exponential tails, that is,
%
%
\begin{equation}
\label{exptail} \exists\kappa>0, C<\infty\mbox{ such that } p(x) \le
Ce^{-\kappa|x|}\ \forall x\in{\mathbb{Z}^d}.
\end{equation}\vfill\eject\noindent
Here $|(x_1,\ldots,x_d)|=\max_i|x_i|$. We define the
\emph{local density} $f_i=f_i(x,\xi)$ of type $i$ near $x\in\Z^d$
by
\[
f_i(x,\xi)=\sum_{y\in{\mathbb{Z}^d}}p(y-x)1 \bigl\{
\xi(y)=i \bigr\},\qquad i=0,1.
\]

Given $p(x)$ satisfying \eqref{passump} and nonnegative
parameters $(\alpha_0,\alpha_1)$, the sto\-chastic
Lotka--Volterra model of
\cite{NP}, $\mathrm{LV}(\alpha_0,\alpha_1)$, is the
spin-flip process $\xi_t$ with rate function
$c_{\LV}(x,\xi)$ given by
%
%
\begin{equation}
\label{LVrates} c_{\LV}(x,\xi) = \cases{ f_1(x,\xi)
\bigl(f_0(x,\xi) + \alpha_0f_1(x,\xi)
\bigr),&\quad $\mbox{if }\xi(x)=0$,\vspace*{2pt}
\cr
f_0(x,\xi)
\bigl(f_1(x,\xi) + \alpha_1f_0(x,\xi)
\bigr),&\quad $ \mbox{if }\xi(x)=1.$}
\end{equation}
All the spin-flip rate functions we will consider, including $c_{\LV}$,
will satisfy the hypothesis
of Theorem B.3 in \cite{Lig99}. By that result, for such a rate
function $c(x,\xi)$, there is a unique $\{0,1\}^{\Z^d}$-valued
Feller process $\xi_t$ with generator equal to the closure of
$\Omega f(\xi)=\sum_{x\in\Z^d}c(x,\xi)(f(\xi^x)-f(\xi))$ on
the space of functions $f$ depending on finitely many
coordinates of $\xi$. Here
$\xi^x$ is $\xi$ but with the coordinate at~$x$
flipped.

One goal of \cite{NP} was to establish coexistence for
$\mathrm{LV}(\alpha_0,\alpha_1)$ for some $\alpha_i$. If we let
$|\xi|=\sum_{x\in{\mathbb{Z}^d}}\xi(x)$ and $\hat\xi(x)=1-\xi
(x)$ for
all $x\in{\mathbb{Z}^d}$, then
coexistence for a spin-flip process $\xi_t$ means that there
is a stationary distribution $\mu$ for $\xi_t$ such that
%
%
\begin{equation}
\label{coexist} \mu\bigl(| \xi| = |\hxi|=\infty\bigr) = 1.
\end{equation}
In \cite{NP}, coexistence was proved for
%
%
\begin{equation}
\label{diag} \alpha=\alpha_0=\alpha_1\in[0,1)
\end{equation}
close enough to 0 and $p(x)=1_{\cN}(x)/|\cN|$, where
%
%
\begin{equation}
\label{NPN} \cN=\bigl\{x\in{\mathbb{Z}^d}\dvtx0<|x|\le L\bigr\}, \qquad L
\ge1,
\end{equation}
excluding only the case $d=L=1$.

A special case of $\mathrm{LV}(\alpha_0,\alpha_1)$ is the voter
model. If we set $\alpha_0=\alpha_1=1$ and use $f_0+f_1=1$,
then $c_{\LV}(x,\xi)$ reduces to the rate function of the
voter model,
%
%
\begin{equation}
\label{VMrates} c_{\vm}(x,\xi) = \bigl(1-\xi(x) \bigr)f_1(x,
\xi) + \xi(x) f_0(x,\xi).
\end{equation}
It is well known (see Chapter~V of \cite{Lig},
Theorems~V.1.8 and V.1.9 in particular) that coexistence for
the voter model is dimension dependent. Let $\zero$
(resp., $\one$) be the element of $\{0,1\}^{\mathbb{Z}^d}$ which
is identically 0 (resp., 1), and let $\delta_\zero$,
$\delta_\one$ be the corresponding unit point masses. If
$d\le2$, then there are exactly two extremal stationary
distributions, $\delta_\zero$ and $\delta_\one$, and hence
no coexistence. If $d\ge3$, then there is a one-parameter
family $\{P_u,u\in[0,1]\}$ of translation invariant extremal
stationary distributions, where $P_u$ has density $u$, that is,
$P_u(\xi(x)=1)=u$. For $u\ne0,1$, each $P_u$ satisfies
\eqref{coexist}, so there is coexistence.

Returning to the general Lotka--Volterra model, coexistence
for $\mathrm{LV}(\alpha_0,\alpha_1)$ for certain
$(\alpha_0,\alpha_1)$ near $(1,1)$ (including $\alpha_0=\alpha_1<1$,
$1-\alpha_i$ small enough) was obtained in
\cite{CP07} for $d\ge3$ and in \cite{CMP} for $d=2$. The
methods used in this work require symmetry in the dynamics
between $0$'s and $1$'s, that is, condition \eqref{diag}. Under
this assumption, Theorem~4 of \cite{CP07} and Theorem~1.2 of
\cite{CMP} reduce to the following, with $\mathrm{LV}(\alpha)$
denoting the Lotka--Volterra model when \eqref{diag} holds.

\renewcommand{\thetheorem}{\Alph{theorem}}
\begin{theorem}\label{thmA}
Assume $d\ge2$
and \eqref{diag} holds. If $d=2$, assume also that
$\sum_{x\in\Z^2}|x|^3p(x)<\infty$. Then there exists
$\alpha_c=\alpha_c(d)<1$ such that coexistence holds for
$\mathrm{LV}(\alpha)$ for $\alpha\in(\alpha_c,1)$.
\end{theorem}

Given coexistence, one would like to know if
there is more than one stationary distribution satisfying
the coexistence condition \eqref{coexist}, if so, what are all
stationary distributions and from what initial states is
there weak convergence to a given stationary distribution.
To state our answers to these questions for $\mathrm{LV}(\alpha)$ we
need some additional notation. Define the hitting times
\[
\tau_\zero=\inf\{t\ge0\dvtx\xi_t=\zero\},\qquad
\tau_\one=\inf\{ t\ge0\dvtx\xi_t=\one\},
\]
and the probabilities, for $\xi\in\{0,1\}^{\mathbb{Z}^d}$,
\begin{eqnarray*}
\beta_0(\xi)& =& P_\xi(\tau_{\zero}<\infty),\qquad
\beta_1(\xi) =P_\xi(\tau_{\one}<\infty),\\
\beta_\infty(\xi) &=& P_\xi(\tau_{\zero}=
\tau_{\one}=\infty),
\end{eqnarray*}
where $P_\xi$ is the law of our process starting at $\xi$.
The point masses $\delta_\zero,\delta_\one$ are clearly
stationary distributions for $\mathrm{LV}(\alpha)$. We write $\xi_t\To\mu$
to mean that the law of $\xi_t$ converges weakly to the
probability measure $\mu$. A law $\mu$ on $\{0,1\}^{\Z^d}$
is symmetric if and only if $\mu(\xi\in\cdot)=\mu(\hat\xi\in
\cdot)$.

We note here that for any translation invariant spin-flip
system $\xi_t$ satisfying the hypothesis (B4) of Theorem~B.3 in
\cite{Lig99},
%
\begin{equation}
\label{betanonzero} \beta_0(\xi)=0\qquad\mbox{if }|\xi|=\infty\quad\mbox{and}\quad \beta_1(\xi)=0\qquad\mbox{if } |\hat\xi|=\infty.
\end{equation}
To see this for $\beta_0$, assume $\xi_0$ satisfies $|\xi_0|=\infty$.
By assumption, there is a
uniform maximum flip
rate $M$ at all sites in all
configurations. So for $\xi_0$ and
$x$ such that $\xi_0(x)=1$,
$P(\xi_t(x)=1)\ge e^{-Mt}$. Since $|\xi_0|=\infty$, we may choose
$A_n\subset{\mathbb{Z}^d}$ satisfying
$|A_n|= n$, $\min\{|x-y|\dvtx x,y\in A_n, x\ne y\}\to\infty$
as $n\to\infty$, and $\xi_0(x)=1\
\forall x\in A_n$. Our hypotheses and translation invariance allow us
to apply Theorem~I.4.6 of \cite{Lig} and conclude that for any
fixed $t>0$, $E (\prod_{x\in A_n}\hat\xi_t(x) )
-\prod_{x\in A_n}E(\hat\xi_t(x))\to0$.\vadjust{\goodbreak}
It follows that for any $n$, there are $\{\vep_n\}$ approaching $0$ so that
\begin{eqnarray*}
P \bigl(\xi_t(x)=0\ \forall x\in\xi_0 \bigr) &\le& P
\bigl(\xi_t(x)=0\ \forall x\in A_n \bigr)
\\
&\le&\vep_n+\prod_{x\in A_n}P \bigl(
\xi_t(x)=0 \bigr)
\\
&\le&\vep_n+ \bigl(1-e^{-Mt} \bigr)^n
\rightarrow0\qquad \mbox{as }n\to\infty.
\end{eqnarray*}

Recall (see Corollary~V.1.13 of \cite{Lig}) that for the voter model
itself and $\xi_0$ translation invariant with $P(\xi_0(x)=1)=u$, we
have $\xi_t\To u\delta_\one+(1-u)\delta_\zero$ if $d\le2$, and
$\xi
_t\To P_u$ if $d\ge3$ and $\xi_0$ is ergodic.
%
%
\begin{thmm}\label{thmLVCCT} Assume
$d\ge2$, and \eqref{exptail}. There exists
$\alpha_c<1$ such that for all
$\alpha\in(\alpha_c,1)$, $\mathrm{LV}(\alpha)$ has a unique translation invariant symmetric
stationary distribution $\nu_{1/2}$
satisfying the coexistence property \eqref{coexist}, such
that for all $\xi_0\in
\{0,1\}^{{\mathbb{Z}^d}}$,
%
%
\begin{equation}
\label{eqLVCCT} \xi_t \To\beta_0(\xi_0)
\delta_\zero+ \beta_\infty(\xi_0)
\nu_{1/2} +\beta_1(\xi_0)\delta_{\one}\qquad
\mbox{as } t\to\infty.
\end{equation}
\end{thmm}

Theorem~\ref{thmLVCCT} is a \emph{complete convergence
theorem}, it gives complete answers to the questions
raised above. The first theorem of this type for infinite
particle systems was proved for the contact process in
\cite{G78}, where $\beta_1(\xi)=0$ for $\xi\ne
\one$ and $\delta_\one$ is
not a stationary distribution. Our result is closely akin
to the complete convergence theorem proved in \cite{H99} for
the threshold voter model. (Indeed, we make use of a number of
ideas from \cite{H99}.) A more recent example is Theorem 4
in \cite{SS} for the $d=1$ ``rebellious voter model.'' For
$p(x)=1_{\cN}(x)/|\cN|$, $\cN$ as in \eqref{NPN}, the
existence and uniqueness of $\nu_{1/2}$ in the above context
follows from results in~\cite{SS} and Theorem \ref{thmA}. The
relationship between Theorem~\ref{thmLVCCT} and results in
\cite{SS} is discussed further in Remarks~\ref{LVAV} and \ref{SScomp}
below.

For $\mathrm{LV}(\alpha)$, we
note that if $0<|\xi_0|<\infty$, then $0<\beta_0(\xi_0)<1$,
where the upper bound is valid for $\alpha$ close enough
to $1$, and if $|\xi_0|=\infty$, then $\beta_0(\xi_0)=0$. By
the symmetry condition \eqref{diag}, this implies that the
obvious symmetric statements with $(\hxi_0,\beta_1)$ in place
of $(\xi_0,\beta_0)$ also hold by \eqref{diag}.
To see the above, note first that
$|\xi_0|<\infty$ trivially implies $\beta_0>0$ since one can
prescribe a finite sequence of flips that leads to the trap
$\zero$. The fact that $\beta_0<1$ for $\alpha<1$ close
enough to $1$ follows from the survival results in
\cite{CP07} for $d\ge3$ (see Theorem~1 there), and in
\cite{CMP} for $d=2$ (see Theorem~1.4 there). Finally, $\beta_0(\xi
_0)=0$ if $|\xi_0|=\infty$ holds by \eqref{betanonzero}.

As our earlier comments on the ergodic theory of the voter
model show, the situation is quite different for $\alpha=1$
as \eqref{eqLVCCT} does not hold. Moreover, by constructing
blocks of alternating $0$'s and $1$'s on larger and larger
annuli one can construct an initial $\xi_0\in\{0,1\}^{\mathbb{Z}^d}$
for which the law of $\xi_t$ does not converge as
$t\to\infty$. This suggests that the above theorem is rather
delicate. Nonetheless we make the following conjecture:

%
\begin{conj} For $\alpha_i<1$, close enough to $1$ and with $\alpha
=(\alpha_0,\alpha_1)$ in the coexistence region of Theorem~1.10 of
\cite
{CDP11}, the
complete convergence theorem holds with a unique nontrivial
stationary distribution $\nu_\alpha$ in place of~$\nu_{1/2}$.
\end{conj}

If $\alpha$ approaches $(1,1)$
so that $\frac{1-\alpha_1}{1-\alpha_0}\to m$, then by Theorem~1.10 of
\cite{CDP11}, the limiting particle density of $\nu_\alpha$ must
approach $u^*(m)$, where $u^*$ is as in (1.50) of~\cite{CDP11}. Hence
one obtains the one-parameter family of invariant laws for the voter
model in the limit along different slopes approaching $(1,1)$.

The $d\ge3$ case of Theorem~\ref{thmLVCCT} is a special
case of a general result for certain \emph{voter
model perturbations}. We will define this class following
the formulation in~\cite{CDP11} (instead of \cite{CP05}),
and then give the additional required definitions needed for
our general result. A \emph{voter model perturbation} is a
family of spin-flip systems
$\xi^\vep_t$, $0<\vep\le\vep_0$ for some $\vep_0>0$, with
rate functions
%
%
\begin{equation}
\label{vmpert} c_\vep(x,\xi) = c_{\vm}(x,\xi) +
\vep^2 c^*_\vep(x,\xi) \ge0,\qquad x\in\Z^d, \xi\in
\{0,1\}^{\Z^d},
\end{equation}
where $c^*_\vep(x,\xi)$ is a translation invariant, signed
perturbation of the form
%
%
\begin{equation}
\label{vmpert2} c^*_\vep(x,\xi)= \bigl(1-\xi(x) \bigr)h_1^\vep(x,
\xi) + \xi(x)h_0^\vep(x,\xi).
\end{equation}
Here we assume \eqref{passump} and \eqref{exptail} hold, and
for some finite $N_0$ there is a law $q_Z$ of
$(Z^1,\ldots,Z^{N_0})\in\Z^{dN_0}$, functions $g^\vep_i$ on
$\{0,1\}^{N_0}$, $i=0,1$ and $\vep_1\in(0,\infty]$
so that $g_i^\vep\ge0$,
and for $i=0,1$, $\xi\in\{0,1\}^{\Z^d}$, $x\in\Z^d$ and $\vep\in
(0,\vep_1)$,
%
%
\begin{equation}
\label{hgrepn}h_i^\vep(x,\xi) = -\vep_1^{-2}f_i(x,
\xi)+E_Z \bigl(g_i^\vep\bigl(\xi\bigl(x+
Z^1 \bigr), \ldots,\xi\bigl(x+ Z^{N_0} \bigr) \bigr) \bigr).
\end{equation}
Here $E_Z$ is
expectation with respect to $q_Z$.
We also suppose that
(decrease $\kappa>0$ if necessary)
%
%
\begin{equation}
\label{expbd2} P \bigl( Z^*\ge x \bigr) \le C e^{-\kappa x}\qquad\mbox{for }x>0,
\end{equation}
where $Z^* = \max\{ |Z^1|, \ldots,|Z^{N_0}| \}$,
and
there are limiting maps $g_i\dvtx\{0,1\}^{N_0}\to\R_+$ such that for some
$c_g,r_0>0$,
%
%
\begin{equation}
\label{gcvgce} \bigl\Vert g_i^\vep-g_i
\bigr\Vert_\infty\le c_g\vep^{r_0},\qquad i=0,1.
\end{equation}
In addition, we will always assume that for $0<\vep\le\vep_0$,
%
%
\begin{equation}
\label{traps} \zero\mbox{ is a trap for }\xi^\vep_t,
\mbox{ that is, }c_\vep(x,\zero)=0.
\end{equation}
In adding \eqref{gcvgce} and \eqref{traps} to the definition of
\textit{voter model perturbation} we have taken some liberty with the
definition in \cite{CDP11}, but these conditions do appear later in
that work for all the results to hold.

It is easy to check that $\mathrm{LV}(\alpha_0,\alpha_1)$ is a voter model
perturbation, as is done in Section~1.3 of~\cite{CDP11}. We
will just note here that if
$\alpha_i=\alpha_i^\vep=1+\vep^2\theta_i $,
$\theta_i\in\R$ and
$h^\vep_i(x,\xi)=\theta_{1-i}f_i(x,\xi)^2$, $i=0,1$, then
$c_{\LV}(x,\xi)$ has the form given
in \eqref{vmpert} and~\eqref{vmpert2}.
Additional examples of voter model perturbations are given
in Section~1 of \cite{CDP11}.\vadjust{\goodbreak} In fact, many interesting models from the
life sciences and social sciences reduce to the voter model
for a specific choice of parameters, and thus in many cases
can be viewed as voter model perturbations.

Coexistence results for voter model perturbations are given
in \cite{CDP11} and \cite{CP07} for $d\ge3$ (and for the two-dimensional
Lotka--Volterra model in \cite{CMP}). Here we will
additionally require that our voter model perturbations be
\emph{cancellative processes}, which we now define
following Section~III.4 of \cite{Lig}; see also Chapter~III
of \cite{G79}. Let $Y$ be the collection of finite subsets
of ${\mathbb{Z}^d}$ and for $x\in{\mathbb{Z}^d}$, $\xi\in\{0,1\}
^{\mathbb{Z}^d}$ and $A\in
Y$, let $H(\xi,A)=\prod_{a\in A}(2\xi(a)-1)$ (an empty
product is 1). We will call a translation invariant flip
rate function $c(x,\xi)$ (not necessarily a voter model
perturbation) cancellative if there is a positive constant
$k_0$ and a map $q_0\dvtx Y\to[0,1]$ such that
%
%
\begin{equation}
\label{q00} c(x,\xi) = \frac{k_0}{2} \biggl( 1- \bigl(2\xi(x)-1 \bigr)
\sum
_{A\in Y}q_0(A-x)H(\xi,A) \biggr),
\end{equation}
where $A-x=\{a-x\dvtx a\in A\}$, $q_0(\varnothing)=0$,
%
%
\begin{eqnarray}
\sum_{A\in Y} q_0(A)&=& 1\quad \mbox{and }
\label{q01}
\\
\sum_{A\in Y} |A|q_0(A)&<&\infty.
\label{q02}
\end{eqnarray}
This is a subclass of the corresponding processes defined in
\cite{Lig}. It follows from~\eqref{q01} that $c(x,\one)=0$
and so $\one$ is a trap for $\xi$.
The above rate will satisfy the hypothesis
of Theorem B.3 in \cite{Lig99} and so, as discussed above, determines
a unique $\{0,1\}^{\Z^d}$-valued
Feller process; see the discussion in Section III.4 of \cite{Lig}
leading to (4.8) there. (One can also check easily that the same is
true of
our voter model perturbations but at times we will only assume the
above cancellative property.)

Given $c(x,\xi), k_0,q_0$
as above, we can define a
continuous time Markov chain taking values in $Y$ by the following.
For $F,G\in Y$, $F\neq G$, define
%
%
\begin{equation}
\label{qmatrix} Q(F,G) = k_0\sum_{x\in F}
\sum_{A\in Y} q_0(A-x) 1 \bigl\{ \bigl(F
\setminus\{x\} \bigr)\Delta A=G \bigr\},
\end{equation}
where $\Delta$ is the symmetric
difference operator. As noted in \cite{Lig}, $Q$ is the
$Q$-matrix of a nonexplosive Markov chain $\zeta_t$ taking
values in $Y$; see also \cite{G79}. If we think of~$\zeta_t$ as the set of sites occupied by a system of
particles at time $t$, then the interpretation of
\eqref{qmatrix} is this. If the current state of the chain
is $F$, then at rate $k_0$ for each $x\in F$:
\begin{longlist}[(1)]
\item[(1)] $x$ is removed from $F$, and
\item[(2)] with probability $q_0(A-x)$, particles are sent
from $x$ to $A$, with the proviso that a particle landing
on an occupied site $y$ annihilates itself and the
particle at $y$.\vadjust{\goodbreak}
\end{longlist}
Perhaps the simplest example of a cancellative/annihilative
pair $(\xi_t,\zeta_t$) is the voter model and its dual annihilating
random walk system. Here $c_{\vm}(x,\xi)$ satisfies~\eqref{q00}
with $k_0=1$, $q_0(\{y\})=p(y)$, $q_0(A)=0$ if $|A|>1$; again, see
\cite
{G79} and \cite{Lig}. A second example, as
shown in \cite{NP}, is the Lotka--Volterra process, assuming
\eqref{diag} and $p(x)=1_\cN(x)/|\cN|$, $\cN$ satisfies
\eqref{NPN} [this will be extended to our general $p(\cdot)$'s in
Section~\ref{secalmostlast}].

The Markov chain $\zeta_t$ is the \emph{annihilating dual}
of $\xi_t$. The general duality equation of
Theorem~III.4.13 of \cite{Lig} (see also Theorem~III.1.5 of
\cite{G79}) and \cite{Lig}, simplifies in the current setting
to the following \textit{annihilating duality} equation:
%
%
\begin{equation}
\label{eqduality0} E \bigl(H(\xi_t,\zeta_0) \bigr) =
E \bigl(H(\xi_0,\zeta_t) \bigr)\qquad \forall
\xi_0\in\{0,1\}^{\mathbb
{Z}^d}, \zeta_0\in Y.
\end{equation}
In Section~\ref{secH} we will recall from \cite{G79} and
\cite{Lig} several implications of this duality equation
for the ergodic theory of $\xi_t$.

Let $Y_e$ (resp., $Y_o$) denote the set of
$A\in Y$ with $|A|$ even
(resp., odd). We call $\zeta_t$ (or $Q$) \emph{parity
preserving} if
%
%
\begin{equation}
\label{parity} Q(F,G)=0 \qquad\mbox{unless }F,G\in Y_e \mbox{ or }
F,G\in Y_o.
\end{equation}
Clearly $\zeta_t$
is parity preserving if and only if $q_0(A)=0$ for all $A\in Y_e$.
If $\zeta_t$ is parity preserving we will call $\zeta_t$
\emph{irreducible} if $\zeta_t$ is irreducible on $Y_o$ and
also on $Y_e\setminus\{\varnothing\}$, and $Q(A,\varnothing)>0$
for some $A\ne\varnothing$.

One fact we need now is
Corollary~{III.1.8} of \cite{G79}. Let $\mu_{1/2}$ be
Bernoulli product measure on $\{0,1\}^{\mathbb{Z}^d}$ with density
$1/2$. Then under \eqref{traps} there is a translation invariant distribution
$\nu_{1/2}$ with density $1/2$ such that
%
%
\begin{equation}
\label{nuhalf1} \mbox{if the law of $\xi_0$ is $
\mu_{1/2}$ then }\xi_t\To\nu_{1/2} \mbox{ as }t
\to\infty;
\end{equation}
see \eqref{zetaconv} below for a proof. For a
cancellative process, $\nu_{1/2}$ will always denote this
measure. We note that $\nu_{1/2}$ might be $\frac12
(\delta_\zero+\delta_\one)$ and hence not have the
coexistence property \eqref{coexist}.

Theorem~1.15 of \cite{CDP11} gives conditions which
guarantee coexistence for $\xi^\vep_t$ for small positive
$\vep$. One assumption of that result, which we will need
here, requires a function $f$ defined in terms of the voter
model equilibria $P_u$ previously introduced. For bounded
functions $g$ on $\{0,1\}^{\mathbb{Z}^d}$ write $\langle g\rangle
_u=\int g(\xi)\,dP_u(\xi)$, and note that $\langle
g(\xi)\rangle_u =\langle g(\hxi)\rangle_{1-u}$. As in
\cite{CDP11}, define
%
%
\begin{equation}
\label{f} f(u) = \bigl\langle\bigl(1-\xi(0) \bigr)c^*(0,\xi) - \xi
(0)c^*(0,\xi)
\bigr\rangle_u,\qquad u\in[0,1],
\end{equation}
where $c^*$ is as in \eqref{vmpert2} but with $g_i$ in place of
$g_i^\vep$.
As noted in Section 1 of \cite{CDP11}, $f$ is a polynomial
of degree at most $N_0+1$, and is a cubic for
$\mathrm{LV}(\alpha_0,\alpha_1)$.

We extend our earlier definitions of $\beta_i$ and $\tau_i$ to general
spin-flip processes~$\xi$.

\begin{definition*}[(Complete convergence)]
We say
that \textit{the complete convergence theorem} holds for a
given cancellative process $\xi_t$ if \eqref{eqLVCCT} holds
for all initial states $\xi_0\in\{0,1\}^{\Z^d}$, where
$\nu_{1/2}$ is given in \eqref{nuhalf1},
and that it holds \textit{with coexistence} if, in addition, $\nu
_{1/2}$ satisfies
\eqref{coexist}.
\end{definition*}


%
\begin{thmm}\label{thmCCTpert} Assume
$d\ge3$, $c_\vep(x,\xi)$ is a voter model perturbation
satisfying \eqref{exptail},
\eqref{q00}--\eqref{q02} and $f'(0)>0$.
Then there
exists $\vep_1>0$ such that if $0<\vep<\vep_1$ the
complete convergence theorem with coexistence holds for
$\xi^\vep_t$.
\end{thmm}

%
\begin{rem}\label{remexpvsattractive} As can be seen in
our proof of Theorem~\ref{thmCCTpert},
it is possible to drop the exponential tail condition
\eqref{exptail} if the voter model perturbations are
attractive, as is the case for $\mathrm{LV}(\alpha)$; see,
for example, (8.5) with $C_{8.3}=1$ in \cite{CP07} for the
latter. To do this one uses the coexistence result in
Section~6 of \cite{CP07} rather than that in Section~6 of
\cite{CDP11}. In particular it follows that in
Theorem~\ref{thmLVCCT} the complete convergence result
holds for the Lotka--Volterra models considered there for
$d\ge3$ without the exponential tail condition
\eqref{exptail}. For $\mathrm{LV}(\alpha)$ with $d=2$ we will have
to use coexistence results in \cite{CMP} to derive the
complete convergence results, and instead of
\eqref{exptail} these results only require
\[
\sum_{x\in\Z^2}|x|^3p(x)<\infty.
\]
See Remark~\ref{noexptail} in Section~\ref{secalmostlast}.


\end{rem}

Theorem~1.3 of \cite{CDP11} states that if the ``initial rescaled
approximate densities of~$1$'s'' approach
a continuous function $v$ in a certain sense, then the rescaled
approximate densities of $\xi_t$ converge to the unique solution of the
reaction diffusion equation
\[
\frac{\partial u}{\partial t}=\frac{\sigma^2}{2}\Delta u+f(u),\qquad u_0=v.
\]
Hence the condition $f'(0)>0$ means there is a positive
drift for the local density of $1$'s when the density of 1's
is very small and so by symmetry a negative drift when the
density of 1's is close to $1$. In this way we see that
this condition promotes coexistence. It also excludes voter
models themselves for which the complete convergence theorem
fails.

We present two additional
applications of Theorem~\ref{thmCCTpert}.

\begin{example}[(Affine voter model)]\label{ex1}
Suppose
%
%
\begin{equation}
\label{N}\qquad  \cN\in Y\mbox{ is nonempty, symmetric and does not contain the
origin.}
\end{equation}
The corresponding \textit{threshold voter model} rate
function, introduced in \cite{CD91}, is
\[
c_{\tvm}(x,\xi) = 1 \bigl\{\xi(x+y)\ne\xi(x)\mbox{ for some } y\in\cN
\bigr\}.
\]
See Chapter~II of \cite{Lig99} for a general treatment of
threshold voter models, and \cite{H99} for
a complete convergence theorem.
The affine voter model with parameter $\alpha\in[0,1]$, $\mathrm{AV}(\alpha)$,
is the spin-flip system with rate function
%
%
\begin{equation}
\label{AVrates} c_{\AV}(x,\xi)=\alpha c_{\vm}(x,\xi)+(1-
\alpha)c_{\tvm}(x,\xi),
\end{equation}
where $c_{\vm}$ is as in \eqref{VMrates}.
This model is studied in \cite{SS} with voter kernel $p(x)=1_{\cN
}(x)/|\cN|$, as an example of a competition model where locally rare
types have a competitive advantage.
\end{example}

%
\begin{thmm}\label{thmtvm} Assume $d\ge3$, \eqref{exptail}
holds and $\cN$ satisfies \eqref{N}. There is an
$\alpha_c\in(0,1)$ so that for all
$\alpha\in(\alpha_c,1)$, the complete convergence theorem
with coexistence holds for $\mathrm{AV}(\alpha)$.
\end{thmm}
%
%
\begin{rem}\label{LVAV}
It was shown in Theorem~3(a) of \cite{SS} that, excluding
the case $d=1$ and $\cN=\{-1,1\}$, if
$p(x)=1_{\cN}(x)/|\cN|$, $\cN$ as in \eqref{NPN}, and coexistence holds
for $\mathrm{LV}(\alpha)$, respectively, $\mathrm{AV}(\alpha)$, for a given $\alpha<1$,
then there is a unique translation invariant stationary
distribution $\nu_{1/2}$ satisfying \eqref{coexist}. Hence this is
true for $\mathrm{LV}(\alpha)$ in $d\ge2$ for $\alpha<1$, and sufficiently
close to $1$, by Theorem~\ref{thmA}, and for $\alpha$ sufficiently small by
\cite
{NP}. It is also true for $\mathrm{AV}(\alpha)$ for $\alpha=0$ by results in
\cite{CD91} and \cite{Lig94}. The same result in \cite{SS} also shows
that if, in addition, the dual satisfies a certain
``nonstability'' condition, then $\xi_t\To\nu_{1/2}$ if
the law of $\xi_0$ is translation invariant and satisfies
\eqref{coexist}. The complete convergence results in Theorems~\ref
{thmLVCCT} and \ref{thmtvm} above (which are special cases of
Theorem~\ref{thmCCTpert} if $d\ge3$) assert a stronger and
unconditional conclusion for both models for $\alpha$ near $1$.
\end{rem}

\begin{example}[(Geometric voter model)]\label{ex2}
Let $\cN$ satisfy \eqref{N}. The geometric voter model with
parameter $\theta\in[0,1]$, $\operatorname{GV}(\theta)$, is the spin-flip
system with rate function
%
%
\begin{equation}
\label{geom1} c_{\gv}(x,\xi) = \frac{1-\theta^{j}}{1-\theta^{|\cN
|}}\qquad \mbox{if } \sum
_{y\in
\cN}1 \bigl\{\xi(x+y)\ne\xi(x) \bigr\}=j,
\end{equation}
where the ratio is interpreted as $j/|\cN|$ if $\theta=1$.
This \emph{geometric} rate function was introduced in \cite{CD91},
where it was shown to be cancellative. As $\theta$ ranges
from 0 to 1 these dynamics range from the threshold voter
model to the voter model. It turns out that the geometric
voter model is a voter model perturbation for $\theta$ near
1, and the following result is another consequence of
Theorem~\ref{thmCCTpert}.
\end{example}

%
\begin{thmm}\label{thmgeom} Assume $d\ge3$ and $\cN$
satisfies \eqref{N}. There is a $\theta_c\in(0,1)$ so that for all
$\theta\in(\theta_c,1)$, the complete convergence theorem with
coexistence holds for $\operatorname{GV}(\theta)$.
\end{thmm}

%
\begin{rem}[(Comparison with \cite{H99} and
\cite{SS})]\label{SScomp} The emphasis in \cite{SS} was on the use of the annihilating
dual to study the invariant laws and the long time behavior
of cancellative systems. A general result (Theorem 6 of
\cite{SS}) gave conditions on the dual to ensure the
existence of a unique translation invariant stationary law
$\nu_{1/2}$ which satisfies the coexistence property
\eqref{coexist} and a stronger local nonsingularity
property. It also gives stronger conditions on the dual
under which $\xi_t\To\nu_{1/2}$ providing the initial law
is translation invariant and satisfies the above local
nonsingularity condition. The general nature of these
interesting results make them potentially useful in a
variety of settings if the hypotheses can be
verified.

In our work we focus on cancellative systems which are also voter perturbations.
A non-annihilating
dual particle system was constructed in \cite{CDP11} to analyze
the latter, and it is by using both dual processes
that we are able
to obtain a complete convergence theorem
in Theorem~\ref{thmCCTpert} for small perturbations and $d\ge3$
($d\ge2$ for LV in Theorem~\ref{thmLVCCT}).

Theorem~1.1 of \cite{H99} gives a complete convergence theorem for
the threshold voter model, the spin-flip system with rate
function $c_{\tvm}$ given in Example~\ref{ex1} above, and
a complete convergence result \textup{is} established in
\cite{SS} for the one-dimensional ``rebellious voter
model'' for a sufficiently small parameter
value. In both of these works, one fundamental step is to show
that the annihilating dual $\zeta_t$ grows when it survives, a result
we will adapt for use here; see Lemma~\ref{lemH} and the discussion
following Remark~\ref{remHaltcond} below.
Both \cite{H99} and \cite{SS} then use special properties
of the particle systems being studied to complete the proof. In
Proposition~\ref{propCCT} below we give general
conditions under which a cancellative spin-flip system
will satisfy a complete convergence theorem with
coexistence. We then verify the required conditions for the voter
model perturbations
arising in Theorems~\ref{thmCCTpert} and~\ref{thmLVCCT}.
\end{rem}

We conclude this section with a ``flow chart'' of the proof of the main
results, including an outline of the paper.
First, the rather natural condition we impose that~$\mathbf{0}$ is a trap
for our cancellative systems $\xi_t$ will imply that $\xi_t$ is in fact
symmetric with respect to interchange of $0$'s and $1$'s; see
Lemma~\ref
{elemann} in Section~\ref{secH}. This helps explain the asymmetric
\textit{looking} condition $f'(0)>0$ in Theorem~\ref{thmCCTpert} and the
restriction of our results to $\mathrm{LV}$ with $\alpha_1=\alpha_2$.
Section~\ref{secH} also reviews the ergodic theory of
cancellative and annihilating systems.

As noted above, the core of our proof, Proposition~\ref{propCCT}, establishes
a complete convergence theorem for cancellative particle systems (where
$\mathbf{0}$ is a trap), assuming three conditions: (i) growth of the dual
system when it survives, that is, \eqref{eqH}, (ii) a condition
\eqref
{flip2} ensuring a large number of $0$--$1$ pairs at locations separated
by a fixed vector $x_0$ for large $t$ with high probability (ruling out
clustering which clearly is an obstruction to any complete convergence
theorem) and (iii) a condition \eqref{oddgoal} which says if the
initial condition $\xi_0$ contains a large number of $0$--$1$ pairs with
$1$'s in a set $A$, then at time $1$ the probability of an odd number
of $1$'s in $A$ will be close to $1/2$. With these inputs, the proof of
Propostion~\ref{propCCT} in Section~\ref{secflip} is a reasonably
straightforward duality argument. This result requires no voter
perturbation assumptions and may therefore have wider applicability.

We then verify the three conditions for voter model perturbations. The
dual growth condition \eqref{eqH} is established in Lemma~\ref{lemH}
and Remark~\ref{remHaltcond} in Section~\ref{secH}, assuming the dual
is irreducible and the cancellative system itself satisfies $\limsup
_{t\to\infty} P(\xi_t(0)=1)>0$ when $\xi_0=\delta_0$. The latter
condition will be an easy by-product of our percolation arguments in
Section~\ref{secthmproof}. The irreducibility of the annihilating dual
is proved for cancellative systems which are voter model perturbations
in Section~\ref{secirred}; see Corollary~\ref{vpirred}. Condition
\eqref{oddgoal} is verified for voter model perturbations in
Lemma~\ref
{lemodd2} of Section~\ref{secflip}, following ideas in \cite{BDD}. In
Section~\ref{secthmproof} (see Lemma~\ref{lemflip2}) condition
\eqref{flip2} is derived for the voter model perturbations in
Theorem~\ref{thmCCTpert} using a comparison with oriented percolation
which in turn relies on input from \cite{CDP11} (see Lemma~\ref
{vmpertspercolate}) and our condition $f'(0)>0$. Another key in this
argument is the use of certain irreducibility properties of voter
perturbations to help set up the appropriate block events. More
specifically, with positive probability it allows us to transform a
$0$--$1$ pair at a couple of input sites at one time
into a mixed configuration which has a ``positive density'' of both
$0$'s and $1$'s at a later time; see Lemma~\ref{vmirred}. The
percolation comparison will provide a large number of the inputs, and
the mixed configuration will be chosen to ensure a $0$--$1$ pair at
sites with the prescribed separation by $x_0$.
In Section~\ref{secthmproof} we finally prove
Theorem~\ref{thmCCTpert}. Theorem~\ref{thmLVCCT} is
proved in Section~\ref{secalmostlast}, and the proofs of
Theorems~\ref{thmtvm} and \ref{thmgeom} are given in
Section~\ref{seclast}. All of these latter results are proved as
corollaries to Theorem~\ref{thmCCTpert}, except for the
two-dimensional case of Theorem~\ref{thmLVCCT}, where the input for
the percolation argument is derived from \cite{CP07} instead of \cite{CDP11}.


\section{Cancellative and annihilating processes: Growth of the
annihilating dual}\label{secH}
Our main objective in this section (Lemma~\ref{lemH} below) is to show
the dual growth condition: under appropriate hypotheses, the
annihilating dual process $\zeta_t$ will either die out or grow without
bound as $t\to\infty$.

We begin by pointing out the consequences of the assumption that
$\zero$ is a trap for $\xi_t$. We assume here that
$c(x,\xi)$ is a translation invariant
cancellative flip rate function satisfying
\eqref{q00}--\eqref{q02}, $\xi_t$ is the corresponding
cancellative process and $\zeta_t$ the corresponding
annihilating process [the Markov chain on
$Y$ with $Q$-matrix defined in
\eqref{qmatrix}]. In part (iv) below we identify $\xi_t$
with the set of sites of type $1$. Recall that
$H(\xi,A)=\prod_{a\in A}(2\xi(x)-1)$.

%
\begin{lem}\label{elemann}
If $\xi_t$ and $\zeta_t$ are as above, then the following are equivalent:
\begin{longlist}[(iii)]
\item[(i)] $\zero$ is a trap for $\xi_t$.

\item[(ii)] $q_0(A)=0$ for all $A\in Y_e$, that is, $\zeta_t$ is
parity-preserving.\vadjust{\goodbreak}

\item[(iii)] $\xi_t$ is symmetric, that is, $c(x,\xi)=c(x,\hxi)$.

\item[(iv)] The simplified duality equation holds
%
%
\begin{equation}
\label{eqduality2} P\bigl(|\xi_t \cap\zeta_0| \mbox{ is
odd}\bigr) = P\bigl(|\xi_0 \cap\zeta_t| \mbox{ is odd}\bigr)\qquad \forall
\xi_0\in\{ 0,1\}^{\mathbb{Z}^d}, \zeta_0\in Y.
\end{equation}
\end{longlist}
\end{lem}

\begin{pf} Note that $H(\zero, A)=(-1)^{|A|}$, which by
\eqref{q00} implies
\[
c(0,\zero)=\frac{k_0}{2} \biggl(1+\sum_{A\in
Y}q_0(A)
(-1)^{|A|} \biggr).
\]
Thus $\zero$ is a trap for $\xi_t$ if and only if
$ \sum_{A\in Y}q_0(A) (-1)^{|A|}=-1$.
Using \eqref{q01}, we see that
\[
\sum_{A\in Y}q_0(A) (-1)^{|A|}=
\sum_{A\in Y_e}q_0(A) -\sum
_{A\in Y_o}q_0(A) \ge\sum
_{A\in Y_e}q_0(A)-1,
\]
so (i) and (ii) are equivalent.

Using $H(\widehat\xi,A)=(-1)^{|A|}H(\xi,A)$ and
\eqref{q00}, (ii) implies
(iii) because
\begin{eqnarray*}
c(x,\widehat\xi)&=&\frac{k_0}{2} \biggl(1- \bigl(1-2\xi(x) \bigr)\sum
_{A\in
Y_o } q_0(A-x) (-1)^{|A|}H(
\xi,A) \biggr)
\\
&=&\frac{k_0}{2} \biggl(1- \bigl(2\xi(x)-1 \bigr)\sum
_{A\in Y} q_0(A-x)H(\xi,A) \biggr)
\\
&=&c(x,\xi).
\end{eqnarray*}
Conversely, if $c(0,\xi)=c(0,\hxi)$ for all $\xi$, the
previous calculation shows that
\[
\sum_{A\in Y}q_0(A)H(\xi,A) = \sum
_{A\in
Y}q_0(A) (-1)^{|A|+1}H(\xi,A).
\]
Plug in $\xi=\one$ to get
\[
\sum_{A\in Y}q_0(A) = \sum
_{A\in Y_o}q_0(A) - \sum_{A\in Y_e}q_0(A),
\]
which implies $q_0(A)=0$ if $|A|$ is even.
We now have that conditions (i)--(iii) are equivalent.%

The duality equation \eqref{eqduality0} is easily seen to be
equivalent to
\[
P\bigl(|\zeta_0|-|\xi_t\cap\zeta_0|\mbox{ is
odd}\bigr) =P\bigl(|\zeta_t|-|\xi_0\cap\zeta_t|
\mbox{ is odd}\bigr)\qquad \forall\xi_0\in\{0,1\}^{\mathbb{Z}^d},
\zeta_0\in Y.
\]
If $\zeta_t$ is parity preserving, then this is equivalent to
(iv). Conversely, if (iv) holds, and we apply it with
$\xi_0=\zero$ and $\zeta_0=\{x\}$, we get
$P(\xi_t(x)=1)=0$ for
all $t>0$. Since this holds for all $x\in{\mathbb{Z}^d}$,
$\zero$ must be a trap for $\xi_t$.
\end{pf}

We give a brief
review (cf. \cite{G79,Lig}) of the
application of annihilating duality to the ergodic theory of
$\xi_t$. Recall that $\mu_{1/2}$ is Bernoulli product
measure with density $1/2$ on $\{0,1\}^{\mathbb{Z}^d}$.\vadjust{\goodbreak} Let $\zeta^A_t$
denote the Markov chain $\zeta_t$ with initial state $A$,
and let $\xi_0$
have law $\mu_{1/2}$. It is easy to see
by integrating \eqref{eqduality2} with respect to the
law of $\xi_0$ that
%
%
\begin{eqnarray}
\label{zetaconv} P\bigl(|\xi_t\cap A| \mbox{ is odd}\bigr) &=
& E \bigl(P
\bigl(\bigl| \zeta^A_t\cap\xi_0\bigr|\mbox{ is odd }|
\zeta^A_t \bigr)1 \bigl(\zeta^A_t
\neq\varnothing\bigr) \bigr)
\nonumber
\\[-8pt]
\\[-8pt]
\nonumber
&=& \tfrac12 P \bigl(\zeta^A_t\ne\varnothing
\bigr) \qquad\mbox{for all }A\in Y.
\end{eqnarray}
The right-hand side above is monotone in $t$ ($\varnothing$ is
a trap for $\zeta_t$), and so the left-hand side above
converges as $t\to\infty$. By inclusion--exclusion arguments
the class of functions
%
%
\begin{eqnarray}
\label{odddet}\mbox{$ \bigl\{\xi\to1\bigl(|\xi\cap A| \mbox{ is odd}\bigr)\dvtx
A\in Y
\bigr\}$ is a determining class,}
\end{eqnarray}
and hence also a convergence determining class
since the state space is compact.
Therefore the above convergence not only implies \eqref{nuhalf1}, it
characterizes
$\nu_{1/2}$ via: for all
$A\in Y$,
%
%
\begin{equation}
\label{nuhalf2} \nu_{1/2}\bigl(\xi\dvtx|\xi\cap A| \mbox{ is odd}\bigr) =\tfrac12
P \bigl(\zeta^A_t\ne\varnothing\ \forall t\ge0 \bigr).
\end{equation}
The measure $\nu_{1/2}$ is necessarily a translation
invariant stationary distribution for $\xi_t$ with density
$1/2$, and a consequence of \eqref{nuhalf2} is that
$\nu_{1/2}\ne\frac12(\delta_{0} + \delta_{1})$ if and only if
for some $x\ne y\in{\mathbb{Z}^d}$,
%
%
\begin{equation}
\label{dualsurvival} P \bigl(\zeta^{\{x,y\}}_t\ne\varnothing
\ \forall t\ge0 \bigr)>0.
\end{equation}
Thus, a sufficient condition for coexistence for $\xi_t$ is
\eqref{dualsurvival}. Indeed, if \eqref{dualsurvival} holds,
then $\nu_{1/2}(\xi\in\cdot|\xi\notin\{\zero,\one\})$ is a
translation invariant stationary distribution for
$\xi_t$ which must satisfy \eqref{coexist}. (There are
countably many configurations $\xi$ with $|\xi|<\infty$,
none of which can have positive probability because there
are countably many distinct translates of each one.)

Establishing
\eqref{dualsurvival} directly is a difficult problem for
most annihilating systems. [Not so for the annihilating dual
of the voter model, since \eqref{dualsurvival} follows
trivially from transience if $d\ge3$ but fails if $d\le2$.]
To use annihilating duality to go beyond~\eqref{nuhalf1}
requires more information about the behavior of
$\zeta_t$. In particular, one needs that either $|\zeta_t|\to0$ or
$|\zeta_t|\to\infty$ as $t\to\infty$; see \cite{BDD}, for
instance. The following general result gives a condition
for this which we can check for certain voter model
perturbations. It is a key ingredient in the proofs of
Theorems~\ref{thmLVCCT} and~\ref{thmCCTpert}.

We now assume that $\zero$ is a trap for $\xi_t$, and so all the properties
listed in Lemma~\ref{elemann} will hold.

%
\begin{lem}[(Handjani \cite{H99}, Sturm and Swart
\cite{SS})]\label{lemH} Let $\zeta_t$ be
a translation invariant annihilating process
with $Q$-matrix given in \eqref{qmatrix} satisfying
\eqref{q01} and~\eqref{q02}. If $\zeta_t$ is
irreducible, parity-preserving, and satisfies
%
%
\begin{equation}
\label{eqH0} \limsup_{t\to\infty}P \bigl(0\in\zeta_t^{\{0\}}
\bigr) >0,
\end{equation}
then
%
%
\begin{equation}
\label{eqH} \lim_{t\to\infty} P \bigl(0<\bigl|\zeta^B_t\bigr|
\le K \bigr) = 0\qquad  \mbox{for all nonempty }B\in Y\mbox{ and } K\ge1.\vadjust{\goodbreak}
\end{equation}
\end{lem}

%
\begin{rem}\label{remHaltcond} If $\zeta_t$ has
associated cancellative process $\xi_t$ which has $\zero$
as a trap, then the parity-preserving hypothesis in the above
result follows by Lemma~\ref{elemann}. If we let
$\xi^{\{0\}}_t$ denote this process with initial state
$\xi_0^{\{0\}}=\{0\}$, then by the duality
equation \eqref{eqduality2}, \eqref{eqH0} is equivalent
to
%
%
\begin{eqnarray}
\label{eqH3} \limsup_{t\to\infty}P \bigl(\xi^{\{0\}}_t(0)=1
\bigr) > 0.
\end{eqnarray}
\end{rem}

The limit \eqref{eqH} was proved
in \cite{H99} (see Proposition~2.6 there) for the
annihilating dual of the threshold voter model. The
arguments in that work are in fact quite general, and with some work
can be extended to establish Lemma~\ref{lemH} as stated
above. Rather than provide the necessary details, we appeal instead to
Theorem~12 of \cite{SS}, which is proved using a related but
somewhat different approach. To apply this result, and
hence establish Lemma~\ref{lemH}, we must do two
things. The first is to show that (3.54) in \cite{SS} [see \eqref{354}
below] holds;
the second is
to show that our condition \eqref{eqH0} implies the
nonstability condition in Theorem~12 of \cite{SS}. The latter is
nonpositive recurrence
of $\zeta$ ``modulo translations''; see the conclusion of Lemma~\ref
{lemstability} below.

In preparation for these tasks we give a ``graphical
construction'' (as in \cite{G78} or~\cite{CD91}) of $\zeta_t$.
For $x\in\Z^d$, let
$\{(S_n^x,A_n^x,)\dvtx n\in\NN\}$ be the points of independent
Poisson point processes $\{\Gamma^x(ds,dA)\dvtx x\in\Z^d\}$ on $\R
_+\times
Y$ with
rates $k_0\,dsq_0(dA)$. For $R\subset\R^d$ and $0\le t_1\le
t_2$ we let
\[
\cF\bigl(R\times[t_1,t_2] \bigr)=\sigma\bigl(
\Gamma^x |_{\Z^d\times
[t_1,t_2]}\dvtx x\in R \bigr).
\]
Then for $S_i=R_i\times I_i$ as above ($i=1,2$), $\cF(S_1)$
and $\cF(S_2)$ are independent if $S_1\cap S_2=\varnothing$.
At time $S^x_n$ draw arrows from $x$ to $x+y$ for each $y\in
A^x_n\setminus\{0\}$. If $0\notin A^x_n$ put a $\delta$ at
$x$ (at time $S^x_n$). For $x,y\in{\mathbb{Z}^d}$ and $s<t$ we say that
$(x,s)\to(y,t)$ if there is a path from $(x,s)$ to $(y,t)$
that goes across arrows, or up but not through $\delta$'s.
That is, $(x,s)\to(y,t)$ if there are sequences
$x_0=x,x_1,\ldots,x_n=y$ and $s_0=s<s_1<\cdots<s_n
<s_{n+1}=t$ such that:
\begin{longlist}[(ii)]
\item[(i)] for $1\le m\le n$, there is an arrow from
$x_{m-1}$ to $x_{m}$ at time $s_m$;
\item[(ii)] for $1\le m\le n+1$, there are no $\delta$'s in
$\{x_{m-1}\}\times(s_{m-1},s_m)$,
\end{longlist}
and no $\delta$ at $(y,t)$.
For $0\le s<t$, $x,y\in{\mathbb{Z}^d}$ and $B\in Y$ define
\begin{eqnarray*}
N^{(x,s)}_t(y) &=& \mbox{ the number of paths up from }(x,s)
\mbox{ to }(y,t),
\\
\zeta^{B,s}_t &= &\biggl\{y\dvtx\sum
_{x\in B} N^{(x,s)}_t(y) \mbox{ is odd}
\biggr\},
\\
\bar\zeta^{B,s}_t &=& \biggl\{y\dvtx\sum
_{x\in B} N^{(x,s)}_t(y)\ge1 \biggr\},
\end{eqnarray*}
and write $\zeta^B_t$ for $\zeta^{B,0}_t$ and
$\bar\zeta^B_t$ for $\zeta^{B,0}_t$.\vadjust{\goodbreak}

The process $\zeta_t$ is the
annihilating Markov chain on $Y$ with $Q$-matrix as in \eqref{qmatrix}.
The process $\bar\zeta_t$ is additive, meaning
$\bar\zeta^{B,s}_t = \bigcup_{x\in
B}\bar\zeta^{(x,s)}_t
$.
Both $\zeta_t$ and $\bar\zeta_t$ are nonexplosive Markov
chains on $Y$.
Also, it is clear that for every $B\in Y$,
%
%
\begin{equation}
\label{contained} \zeta^{B,s}_t\subset\bar
\zeta^{B,s}_t \qquad \forall0\le s\le t<\infty,
\end{equation}
and also that for any fixed $t>0$,
%
%
\begin{equation}
\label{boundedgrowth} \lim_{K\to\infty}P \bigl(\bar
\zeta^{\{0\}}_u\subset[-K,K]^d\ \forall0\le u \le
t \bigr) = 1.
\end{equation}
Furthermore if $A,B\in Y$ satisfy $\min_{a\in A,b\in B}|a-b|>2K$ and
$t>s\ge0$, then
%
%
\begin{eqnarray}
\label{separated}
\zeta_t^{A\cup B,s}=
\zeta_t^{A,s}\cup\zeta_t^{B,s}
\nonumber
\\[-8pt]
\\[-8pt]
 \eqntext{\mbox{on the event }\bigl\{\bar\zeta^{A,s}_u\subset A+[-K,K]^d,
\bar\zeta^{B,s}_u\subset B+[-K,K]^d \ \forall
s\le u\le t \bigr\}}
\end{eqnarray}
(where $A+B=\{x+y\dvtx x\in A, y\in B\}$).

The following result is key to verifying condition (3.54) of \cite{SS}.

%
\begin{lem}\label{graphrep}
Let $A\in Y$, $r\in\NN$ and $B_m=\{y^m_1,\ldots,y^m_r\}\in Y$ be such
that
$\lim_{m\to\infty} \min_i|y^m_i|=\infty$. If $\zeta^A$ and $\zeta
^{B_m}$ are independent copies of $\zeta$ with the given initial
conditions, then for each $t\ge0$ and $n\in\NN$,
\[
\lim_{m\to\infty}P \bigl(\bigl|\zeta_t^{A\cup B_m}\bigr|=n
\bigr)-P \bigl(\bigl|\zeta_t^A\bigr|+\bigl|\zeta_t^{B_m}\bigr|=n
\bigr)=0.
\]
\end{lem}
\begin{pf} Assume $(\zeta^B_t)$ are constructed as above for $B\in Y$ and
$t\ge0$. For $K\in\NN$ define $\tilde\zeta_t^{B,(K)}$ as $\zeta^B_t$
but now only count paths which are contained in $B+[-K,K]^d$. This
implies that
%
%
\begin{equation}
\label{tildemeas} \tilde\zeta_t^{B,(K)}\mbox{ is }\cF
\bigl( \bigl(B+[-K,K]^d \bigr)\times[0,t] \bigr)\mbox{-measurable}.
\end{equation}
Fix $\vep>0$. By \eqref{boundedgrowth} and the additivity of $\bar
\zeta
_t$ we may choose $K(\vep)\in\NN$ so that if $K\ge K(\vep)$, then
%
%
\begin{eqnarray}
\label{contain}\qquad  && P \bigl(\bar\zeta^A_u\subset
A+[-K,K]^d\mbox{ and }\bar\zeta^{B_m}_u
\subset B_m+[-K,K]^d\mbox{ for all }u\in[0,t] \bigr)
\nonumber
\\[-8pt]
\\[-8pt]
\nonumber
&&\qquad >1-\vep\qquad \mbox{for all }m\in\NN.
\end{eqnarray}
Write $\tilde\zeta_t^B$ for $\tilde\zeta^{B,(K(\vep))}_t$. Choose
$m(\vep)\in\NN$ so that $\min_{a\in A,b\in B_m}|a-b|>2K(\vep)$ for
$m\ge m(\vep)$. It follows from \eqref{separated} and \eqref{contained}
that on the set in \eqref{separated} with $K=K(\vep)$, for $m\ge
m(\vep)$,
%
%
\begin{equation}
\label{zetasum}\bigl|\zeta^{A\cup B_m}_t\bigr|=\bigl|\zeta^A_t\bigr|+\bigl|
\zeta_t^{B_m}\bigr|
\end{equation}
and
%
%
\begin{equation}
\label{tildeeq} \zeta^A_t=\tilde\zeta^A_t
\quad\mbox{and}\quad\zeta_t^{B_m}=\tilde\zeta^{B_m}_t.
\end{equation}
[The latter is an easy check using \eqref{separated}.] We conclude from
the last two results
that
%
%
\begin{equation}
\label{tildesum} P \bigl(\bigl|\zeta^{A\cup B_m}_t\bigr|\neq\bigl|\tilde
\zeta^A_t\bigr|+\bigl|\tilde\zeta^{B_m}_t\bigr|
\bigr)<\vep\qquad\mbox{for }m\ge m(\vep).
\end{equation}

By \eqref{tildemeas} and the choice of $m(\vep)$ we see that $\tilde
\zeta^A_t$ and $\tilde\zeta^{B_m}_t$ are independent for $m\ge
m(\vep
)$. Using this independence and then \eqref{tildesum} we conclude that
\begin{eqnarray*}
\hspace*{-4pt}&&\Biggl|P \bigl(\bigl|\zeta_t^{A\cup B_m}\bigr|=n \bigr)- \Biggl(\sum
_{k=0}^nP \bigl(\bigl|\zeta^A_t\bigr|=k
\bigr)P \bigl(\bigl|\zeta_t^{B_m}\bigr|=n-k \bigr) \Biggr) \Biggr|
\\
\hspace*{-4pt}&&\qquad\le\vep+ \Biggl|P \bigl(\bigl|\tilde\zeta^A_t\bigr|+\bigl|\tilde\zeta
^{B_m}_t\bigr|=n \bigr)- \Biggl(\sum
_{k=0}^nP \bigl(\bigl|\zeta^A_t\bigr|=k
\bigr)P \bigl(|\zeta_t^{B_m}|=n-k \bigr) \Biggr) \Biggr|
\\
\hspace*{-4pt}&&\qquad\le\vep+ \Biggl|\sum_{k=0}^n \bigl[P \bigl(\bigl|
\tilde\zeta^A_t\bigr|=k \bigr)P \bigl(\bigl|\tilde\zeta
^{B_m}_t\bigr|=n-k \bigr)-P \bigl(\bigl|\zeta^A_t\bigr|=k
\bigr)P \bigl(\bigl|\zeta_t^{B_m}\bigr|=n-k \bigr) \bigr] \Biggr|
\\
\hspace*{-4pt}&&\qquad\le3\vep.
\end{eqnarray*}
In the last line we have used \eqref{tildeeq}. The result follows.
\end{pf}

Say that $A,B\in Y_o$ are equivalent if they are translates of
each other, let $\tilde Y_o$ denote the set of equivalence
classes, and (abusing notation slightly) let $\tilde A$
denote the equivalence class containing $A\in Y$. Since the
dynamics of $\zeta$ are translation invariant, for
parity-preserving $\zeta$ we may define $\tilde\zeta_t$ as
the $\tilde Y_o$-valued Markov process obtained by taking
the equivalence class of $\zeta_t$.
The
nonstability requirement of Theorem~12 of \cite{SS} is that
$\tzeta_t$ not be positive recurrent on $\tilde Y_o$.

%
\begin{lem}\label{lemstability} If $\zeta$ is parity-preserving,
irreducible and satisfies
\eqref{eqH0}, then $\tzeta_t$ is not positive recurrent.
\end{lem}

\begin{pf} We use the same arguments as in the proof of
Lemma~2.4 of \cite{H99}.
First, $\zeta_t$ cannot be positive recurrent on $Y_o$. To
check this, we first note that translation invariance
implies
\[
P \bigl(\zeta_t^{\{x\}}=\{x\} \bigr)=P \bigl(
\zeta_t^{\{0\}}=\{0\} \bigr)\qquad \mbox{for all }t\ge0, x\in{
\mathbb{Z}^d}.
\]
If $\zeta_t$ is positive recurrent on $Y_o$, then the limit
$\mu(A)=\lim_{t\to\infty}P(\zeta^B_t=A)$ exists and is
positive for all
$A,B\in Y_o$. Letting $t\to\infty$ above, this implies $\mu(\{ 0\})=
\mu(\{ x\})$ for all $x$ which is impossible, so $\zeta_t$
is not positive recurrent on $Y_o$. A consequence of this is
that for any fixed
$k>0$,
\[
\lim_{t\to\infty} P \bigl(\zeta^{\{0\}}_t
\subset[-k,k]^d \bigr) = 0.
\]

Next, suppose $\tzeta_t$ is positive recurrent on $\tilde
Y_o$, with some stationary
distribution~$\tilde\mu$, which must satisfy
\[
\tilde\mu\bigl(A\in\tilde Y_o\dvtx\operatorname{diam}(A) \le k \bigr)
\to1 \qquad\mbox{as }k\to\infty,\vadjust{\goodbreak}
\]
where $\operatorname{diam}(A)=\max\{|x-y|\dvtx x,y\in A\}$ is well defined for $A\in
\tilde Y_o$.
For any $k,t$, since $\operatorname{diam}(\tzeta^{\{0\}}_t)
=\operatorname{diam}(\zeta^{\{0\}}_t)$, we have
\[
P \bigl(0\in\zeta^{\{0\}}_t \bigr) \le P \bigl(
\zeta^{\{0\}}_t\subset[-k,k]^d \bigr) +P \bigl(
\operatorname{diam} \bigl(\tzeta^{\{0\}}_t \bigr) > k \bigr).
\]
Letting $t\to\infty$ gives
\[
\limsup_{t\to\infty}P \bigl(0\in\zeta^{\{0\}}_t
\bigr) \le\tilde\mu\bigl(A\in\tilde Y_o\dvtx\operatorname{diam}(A)>k
\bigr).
\]
The right-hand side above tends to 0 as $k\to\infty$, so we
have a contradiction to the assumption
\eqref{eqH0}.
\end{pf}

\begin{pf*}{Proof of Lemma~\ref{lemH}} Thanks to the
above lemma we have verified all the hypotheses of
Theorem~12 of \cite{SS} except for their (3.54) which we now
state in our notation:
for each $n\in\Z_+$, $L\ge1$ and $t>0$,
%
%
\begin{eqnarray}
\label{354}\qquad  &&\inf\bigl\{P \bigl(\bigl|\zeta_t^A\bigr|=n \bigr)
\dvtx|A|=n+2 \mbox{ and } 0<|i-j|\le L \mbox{ for some }i,j\in A \bigr
\}
\nonumber
\\[-8pt]
\\[-8pt]
\nonumber
&&\qquad>0.
\end{eqnarray}
Assume \eqref{354} fails. Then for some $n$, $L$ and $t$
as above, by translation invariance and compactness of $Y$
(with the subspace topology it inherits from
$\{0,1\}^{\Z^d}$), there are $\{A_m\}\subset Y$ so that for
some integer $2\le s\le n+2$ and $x_2\in[-L,L]^d$,
$A_m=\{0,x_2,\ldots,x_s\}\cup\{x_{s+1}^m,\ldots,x^m_{n+2}\}\equiv
A\cup B_m$, where $\lim_{m\to\infty}|x_i^m|=\infty$ for each
$i\in\{s+1,\ldots,n+2\}$ and
%
%
\begin{equation}
\label{fail354}\lim_{m\to\infty}P \bigl(\bigl|\xi_t^{A_m}\bigr|=n
\bigr)=0.
\end{equation}
By the irreducibility of $\zeta$, $P(|\zeta_t^{A}|=s-2)=p>0$.
If $\zeta_t^{A}$ and $\zeta_t^{B_m}$ are as in Lemma~\ref{graphrep},
then by
that result,
\begin{eqnarray*}
\lim_{m\to\infty}P \bigl(\bigl|\zeta_t^{A_m}\bigr|=n
\bigr)&=&\lim_{m\to\infty}P \bigl(\bigl|\zeta^A_t\bigr|+\bigl|
\zeta^{B_m}_t\bigr|=n \bigr)
\\
&\ge &P \bigl(\bigl|\zeta^A_t\bigr|=s-2 \bigr)\liminf
_{m\to\infty}P \bigl(\bigl|\zeta_t^{B_m}\bigr|=|B_m|
\bigr)
\\
&\ge& p\exp\bigl\{-k_0(n+2-s)T \bigr\}>0.
\end{eqnarray*}
In the last line we use the fact that by its graphical construction,
$\zeta^{B_m}$ will
remain constant up to time $t$ if none of the $|B_m|$ independent rate
$k_0$ Poisson processes attached to each of the sites in $B_m$ fire by
time $t$. This contradicts \eqref{fail354} and so~\eqref{354} must
hold. We now may apply Theorem~12 of \cite{SS} to obtain the required
conclusion.
\end{pf*}

\section{Irreducibility}\label{secirred}
In addition to the explicit irreducibility requirement for $\zeta_t$ in
Lemma~\ref{lemH}, some arguments in
Section~\ref{secthmproof} will require
irreducibility type conditions for the voter model
perturbations $\xi^\vep_t$. We collect and prove the necessary
results for both processes in this section.\vadjust{\goodbreak}

Assuming $\sum_{y\in{\mathbb{Z}^d}} q_0(\{y\})>0$, define the step
distribution of a random walk associated with $q_0$
by
\[
q(x)=q_0 \bigl(\{x\} \bigr) \Big/\sum_{y\in{\mathbb{Z}^d}}
q_0 \bigl(\{y\} \bigr).
\]

%
\begin{lem}\label{suffindecomp} Let $\zeta_t$ be a
parity-preserving annihilating process with $Q$-matrix
given in \eqref{qmatrix}. Assume $q_0(A_0)>0$ for some
$A_0\in Y$ with $|A_0|\ge3$, and for some symmetric,
irreducible random walk kernel $r$ on ${\mathbb{Z}^d}$, $q(x)>0$
whenever $r(x)>0$. Then $\zeta_t$ is irreducible.
\end{lem}

\begin{pf}
The proof is elementary but awkward, so we will only
sketch the argument. Note that if $x\in A$ and $y\notin A$,
then
\[
Q \bigl(A, \bigl(A\setminus\{x\} \bigr)\cup\{y\} \bigr)\ge q_0
\bigl( \{y-x\} \bigr)=cq(y-x).
\]
So by using only the $q_0(\{x\})$
``clocks'' with $r(x)>0$, $\zeta_t$ can with positive
probability execute exactly any finite sequence of
transitions that the annihilating random walk system with
step distribution $r$ can. We will refer to ``$r$-random walks''
below in describing such transitions.

We first check that the assumptions on $q_0$ imply that
$\zeta_t$ can reach any set $B$ with $|B|=|\zeta_0|$ with
positive probability. To see this, we first construct a
set $B'$ by starting $r$-random walks at each site of $B$ and
then moving them apart, one at a time, avoiding
collisions, to widely separated locations, resulting in
$B'$. Note that by reversing this entire sequence of
steps, it is possible to move $r$-random walks starting at the
sites of $B'$ to $B$ without collisions. This uses the
symmetry of $r$. Now, to
move $r$-walks from $\zeta_0$ to $B$ we first move walks from
$\zeta_0$ to some $\zeta_0'$, avoiding collisions, where
the sites of $\zeta_0'$ are widely separated. Pair off
points from $\zeta_0'$ and $B'$ and move $r$-walks one at a
time from $\zeta_0'$ to $B'$ without collisions. This is
possible if $\zeta_0'$ and $B'$ are
sufficiently spread out since $r$ is irreducible. Finally, move the
walks from $B'$
to $B$ without collisions as discussed above.

It should be clear that if $\zeta_0\ne\varnothing$, then
$\zeta_t$ can reach a set $B$ such that $|B|=|\zeta_0|-2$,
since this is the case for annihilating
random walks. Finally, if $\zeta_0\ne\varnothing$, then
$\zeta_t$ can reach a set $B$ with $|B|\ge|\zeta_0|+2$.
Choose $x_1$ far from $\zeta_0$ so that $\zeta_0$
and $x_1+A_0$ are disjoint, and such that for some
$x_0\in\zeta_0$, an $r$-walk starting at $x_0$ can reach $x_1$
by a sequence of steps avoiding $\zeta_0$. Now using the ``$A_0$ clock''
at $x_1$ we get a transition from
$(\zeta_0\setminus\{x_0\})\cup\{x_1\}$ to
$\zeta_0'=(\zeta_0\setminus\{x_0\})\Delta(x_1+A_0)$, and
$|\zeta_0'|\ge|\zeta_0|+2$.
\end{pf}

The next result will allow us to apply the above lemma to
voter model perturbations. Recall $p(x)$ satisfies
\eqref{passump}, and $c_{\vm}(x,\xi)$ is the corresponding
voter model flip rate function.

%
\begin{lem}\label{suffirred} There is an
$\vep_2=\vep_2(p(\cdot))>0$ and $R_1=R_1(p(\cdot))$ such
that:
\begin{longlist}[(ii)]
\item[(i)] $p(\cdot| |x|<R_1)$ is irreducible, and\vadjust{\goodbreak}

\item[(ii)] if
$c(x,\xi)=c_{\vm}(x,\xi)+\tilde c(x,\xi)$ is a
translation invariant, cancellative flip rate function
with $\zero$ as a trap such that
%
%
\begin{equation}
\label{c*cond}\Vert\tilde c\Vert_\infty<\vep_2,\qquad \sum
_{x\neq
0}\bigl|\tilde c(0,\delta_x)\bigr|<
\vep_2,
\end{equation}
then the dual kernel $q_0$ satisfies
%
%
\begin{equation}
\label{qdomp} q_0 \bigl(\{x\} \bigr)> (k_0
3)^{-1}p(x) \qquad\mbox{for all }0<|x|<R_1.
\end{equation}
\end{longlist}
\end{lem}

\begin{pf}
Since $p$ is irreducible, we may choose $R_1$ so that $p(\cdot|
|x|<R_1)$ is also
irreducible. Assume \eqref{c*cond} holds for an appropriate
$\vep_2$ which will be chosen below. We will write $\hat
\xi(B)$ for $\sum_{x\in B}\hat\xi(x)$. Also, \emph{in this proof
only}, we will let $A$ denote a \emph{random set} with probability
mass function
$q_0$, and write $E_0(g(A))=\int g \,dP_0$ for
$\sum_{B\in Y}g(B)q_0(B)$. With this notation, by our hypotheses we have
%
%
\begin{equation}
\label{tworep} c_{\vm}(0,\xi)+\tilde c(0,\xi)=\frac{k_0}{2}
\bigl[1+(-1)^{\xi
(0)}E_0 \bigl((-1)^{\hat\xi(A)} \bigr)
\bigr].
\end{equation}

Recall by Lemma~\ref{elemann} that $P_0(|A|\mbox{ is
odd})=1$. Therefore if we set $\xi=\delta_x$ for $x\neq0$
in \eqref{tworep}, we get
\[
p(x)+\tilde c(0,\delta_x)=\frac{k_0}{2} \bigl[1+E_0
\bigl((-1)^{|A\setminus\{x\}|} \bigr) \bigr]=k_0P_0(x\in A),
\]
and so
%
%
\begin{equation}
\label{P0x}P_0(x\in A)= \bigl(p(x)+\tilde c(0,\delta_x)
\bigr)k_0^{-1}.
\end{equation}

If we take $\xi=\delta_{\{x_0,x_1\}}$ in \eqref{tworep},
where $x_0$, $x_1$ are two distinct nonzero points, then we
get
\begin{eqnarray*}
p(x_0)+p(x_1)+\tilde c(0,\delta_{\{x_0,x_1\}})& =&
\frac{k_0}{2} \bigl[1+E_0 \bigl((-1)^{|A\setminus\{x_0,x_1\}|} \bigr)
\bigr]
\\
&=&k_0P \bigl(1_A(x_0)\neq1_A(x_1)
\bigr),
\end{eqnarray*}
and so
%
%
\begin{equation}
\label{P0x1x2} P_0 \bigl(1_A(x_0)
\neq1_A(x_1) \bigr)= \bigl(p(x_0)+p(x_1)+
\tilde c(0,\delta_{\{
x_0,x_1\}}) \bigr)k_0^{-1}.
\end{equation}

For any two distinct nonzero points, $x_0$ and $x_1$, we have
\[
P_0 \bigl(1_A(x_0)\ne1_A(x_1)
\bigr) = P_0(x_0\in A) + P_0(x_1
\in A) -2 P_0(x_0\in A,x_1\in A).
\]
Therefore, we see that \eqref{P0x} and \eqref{P0x1x2} imply
\begin{eqnarray*}
P_0 \bigl(\{x_0,x_1\}\subset A \bigr) &=&
\tfrac{1}{2} \bigl[P_0(x_0\in A)+P_0(x_1
\in A)-P_0 \bigl(1_A(x_0)\neq
1_A(x_1) \bigr) \bigr]
\\
&=& \bigl[\tilde c(0,\delta_{x_0})+\tilde c(0,\delta_{x_1})-
\tilde c(0,\delta_{\{x_0,x_1\}}) \bigr](2k_0)^{-1},
\end{eqnarray*}
which gives the simple bound
\[
P_0 \bigl(\{x_0,x_1\}\subset A \bigr)\le
\tfrac{3}{2}\|\tilde c\|_\infty k_0^{-1}.
\]

Note that if $0\neq x\in A$ but $A\neq\{x\}$, then $P_0$-a.s. $A$ must
contain $x$ and another nonzero point
as $|A|$ is a.s. odd, and so for $0<|x|<R_1$ and $R_2>R_1$,
\begin{eqnarray*}
&&P_0 \bigl(A=\{x\} \bigr)\\
&&\qquad\ge  P_0(x\in A)-\sum
_{x_1\notin\{0, x\}}P_0 \bigl(\{x,x_1\} \subset A
\bigr)
\\
&&\qquad\ge k_0^{-1} \biggl[p(x)+\tilde c(0,
\delta_x)-(3/2)\Vert\tilde c\Vert_\infty(2R_2+1)^d-
\sum_{|x_1|>R_2}P_0(x_1\in A)
\biggr]
\\
&&\qquad\ge k_0^{-1} \biggl[p(x)- \bigl(1+2(2R_2+1)^d
\bigr)\Vert\tilde c\Vert_\infty-\sum_{|x_1|>R_2}
\bigl(p(x_1)+\tilde c(0,\delta_{x_1}) \bigr) \biggr].
\end{eqnarray*}
We have used the previous displays and \eqref{P0x} in the
above. Recalling the bounds in our assumption \eqref{c*cond} on
$\tilde c$,
we conclude that
%
%
\begin{equation}
\label{qlb} P_0 \bigl(A=\{x\} \bigr)\ge k_0^{-1}
\biggl[p(x)-\sum_{|x_1|>R_2}p(x_1)-2
\bigl(1+(2R_2+1)^d \bigr)\vep_2 \biggr].
\end{equation}
Now let
\[
\eta=\eta\bigl(p(\cdot) \bigr)=\inf\bigl\{p(x)\dvtx|x|<R_1\mbox{
and }p(x)>0 \bigr\}>0,
\]
choose $R_2=R_2(p(\cdot))>R_1$ so that $\sum_{|x_1|>R_2}p(x_1)<\eta/3$
and define
\[
\vep_2=\frac{\eta}{6((2R_1+1)^d+1)}.
\]
Then by \eqref{qlb}
\[
P_0 \bigl(A=\{x\} \bigr)\ge(3k_0)^{-1}p(x)
\qquad\mbox{for all }0<|x|<R_1,
\]
and we are done.
\end{pf}

For the rest of this section we assume $\{\xi^\vep\dvtx0<\vep\le
\vep
_0\}$
is a voter model perturbation
with rate function $c_\vep$ [so that \eqref{vmpert}--\eqref{traps}
are valid] which is also cancellative for each $\vep$ as above with dual
kernels $q_0^\vep$ satisfying \eqref{q00}--\eqref{q02}. In particular
the $\tilde c$ in Lemma~\ref{suffirred} is now
$\vep^2c^*_\vep$. By Lemma~\ref{elemann}, all the
conclusions of that result hold.

%
\begin{cor}\label{vpirred} Assume that
%
%
\begin{equation}
\label{create} \mbox{for small enough $\vep$, } q_0^\vep(A)>0
\mbox{ for some $A\in Y$ with }|A|>1.
\end{equation}
Then there is an $\vep_3>0$ depending on $p$, $\vep_1$, $\{g_i^\vep\}$
and the $\vep$ required in \eqref{create} so that if $0<\vep<\vep_3$,
then the annihilating dual with kernel $q_0^\vep$ is irreducible.
\end{cor}
\begin{pf} Let $R_1$ be as in Lemma~\ref{suffirred}. An easy calculation
shows that
\[
\Vert\tilde c\Vert_\infty\vee\biggl(\sum_x\bigl|
\tilde c(0,\delta_x)\bigr| \biggr)\le\vep^2 \Biggl[
\vep_1^{-2}+\bigvee_{i=0}^1\bigl\Vert
g_i^\vep\bigr\Vert_\infty\Biggr]\le
\vep^2C
\]
for some constant $C$, independent of $\vep$. Therefore for $\vep
<\vep
_3$ ($\vep_3$ as claimed) we have the hypotheses, and hence conclusion,
of Lemma~\ref{suffirred}. This allows us to apply Lemma~\ref
{suffindecomp} with $r(\cdot)=p(\cdot| |x|<R_1)$ and hence conclude
that the annihilating dual $\zeta$ is irreducible for such $\vep$.
\end{pf}
%
%
\begin{rem}\label{nvm} Clearly \eqref{create} is a
necessary condition for the conclusion to hold. In fact
if it fails, it is easy to check that $c_\vep(x,\xi)$ is a
multiple of the voter model rates with random walk kernel
$q^\vep_0(\{x\})$. Hence this condition just eliminates voter
models for which the conclusions of Corollary~\ref{vpirred}, as well as
Lemma~\ref{lemH} and Proposition~\ref{propCCT} below,
will also fail in general.

Note that if \eqref{create} fails, then for some $\vep_n\downarrow0$,
\[
c^*_{\vep_n}(0,\xi)=\vep_n^{-2}c_{\vep_n}(0,
\xi)-\vep_n^{-2}c_{\vm
}(0,\xi)=
\lambda_n\tilde c^n_{\vm}(0,\xi)-
\vep_n^{-2}c_{\vm
}(0,\xi),
\]
where $\tilde c^n_{\vm}(0,\xi)$ is the rate function for the voter
model with kernel $q_0^{\vep_n}(\{\cdot\})$.
From this it is easy to check that if $\langle\cdot\rangle_u$ is
expectation with respect to the voter model equilibrium for $c_{\vm}$
with density $u$, then
\[
\bigl\langle(1-\xi)c^*_{\vep_n}(0,\xi)-\xi c^*_{\vep_n}(0,\xi)
\bigr\rangle_u=0,
\]
and so by \eqref{gcvgce} and \eqref{f}, $f(u)\equiv0$. Therefore,
the condition
$f'(0)>0$ in Theorem~\ref{thmCCTpert} implies
\eqref{create}.
\end{rem}

Next we prove an irreducibility property for the voter model
perturbations $\xi^\vep_t$ themselves. To do so we
introduce the (unscaled) graphical representation for
$\xi^\vep_t$ used in \cite{CDP11}. First put
\[
\bar c= \sup_{\vep<\vep_0} \bigl(\bigl\|g^\vep_1
\bigr\|_\infty+\bigl\|g^\vep_0\bigr\|_\infty+1 \bigr)<
\infty.
\]
For $x\in{\mathbb{Z}^d}$, introduce independent Poisson point
processes on $\R_+$,
$\{
T^{x}_n, n \ge1 \}$ and $\{ T^{*,x}_n, n \ge1\}$, with
rates $1$ and $\vep^2\bar c$, respectively. For $x\in{\mathbb{Z}^d}$
and $n\ge1$, define independent random variables $X_{x,n}$
with distribution $p(\cdot)$, $Z_{x,n}=(Z^1_{x,n}, \ldots,
Z^{N_0}_{x,n})$ with distribution $q_Z(\cdot)$, and $U_{x,n}$
uniform on $(0,1)$. These random variables are independent
of the Poisson processes, and all are independent of any
initial condition $\xi^\vep_0\in\{0,1\}^{\Z^d}$. For
all $x\in{\mathbb{Z}^d}$ we allow $\xi^\vep_t(x)$ to change only at times
$t\in\{T^x_n,T^{*,x}_n,n\ge1\}$. At the voter times $T^x_n,
n \ge1$ we draw a voter arrow from $(x,T^x_n)$ to
$(x+X_{x,n},T^x_n)$ and set
$\xi^\vep_{T^x_n}(x)=\xi^\vep_{T^x_n-}(x+X_{x,n})$. At the times
$T^{*,x}_n$,
$n\ge1$ we draw ``*-arrows'' from $(x,T^{*,x}_n)$\vadjust{\goodbreak}
to each $(x+Z^i_{x,n},T^{*,x}_n)$, $1\le i\le N_0$, and if
$\xi^\vep_{T^{*,x}_n-}(x)=i$ we set $\xi^\vep_{T^{*,x}_N}(x)=1-i$ if
\[
U_{x,n} < g^\vep_{1-i} \bigl(
\xi^\vep_{T^{*,x}_n-} \bigl(x+Z^1_{x,n}
\bigr), \ldots, \xi^\vep_{T^{*,x}_n-} \bigl(x+Z^{N_0}_{x,n}
\bigr) \bigr)/\bar c.
\]
As noted in Section 2 of \cite{CDP11}, this recipe defines a pathwise
unique process $\xi^\vep_t$ whose law is specified by the flip rates
in \eqref{vmpert}. We refer to this as the graphical construction of
$\xi^\vep_t$.
For $x\in\Z^d$, $\{(X_{(x,n)},T_n^x)\dvtx n\in\NN\}$ and $\{
(Z_{x,n},T^{*,x}_n,U_{x,n})\dvtx n\in\NN\}$
are the points of independent collections of independent Poisson
point processes, $(\Lambda^x_w(dy,dt),x\in\Z^d)$ and
$(\Lambda^x_r(dy,dt,du),\break x\in\Z^d)$, on $\Z^d\times\R_+$ with rate
$dt p(\cdot)$, and on $\Z^d\times\R_+\times[0,1]$ with
rate\break  $\vep^2\bar c\,dt q_Z(\cdot) \,du$, respectively. For $R\subset\R^d$
and $0\le t_1\le t_2$
we let
\[
\cG\bigl(R\times[t_1,t_2] \bigr)=\sigma\bigl(
\Lambda^x_w|_{\Z^d\times[t_1,t_2]}, \Lambda^{x'}_r|_{\Z^d\times
[t_1,t_2]\times[0,1]}
\dvtx x,x'\in R \bigr),
\]
that is, the $\sigma$-field generated by the points of the graphical
construction in $R\times[t_1,t_2]$.

A coalescing branching random walk dual for $\xi^\vep_t$ is
constructed in \cite{CDP11}. We give here only the part of
that dual which we need. Using only the Poisson processes
$T^x_n,x\in{\mathbb{Z}^d}$, define a coalescing random walk system as
follows. Fix $t>0$. For each $y\in{\mathbb{Z}^d}$ define $B^{y,t}_u,
u\in[0,t]$ by putting $B^{y,t}_0=y$ and then proceeding
``down'' in the graphical construction and using the voter
arrows to jump. More precisely, if $T^{y}_1>t$ put
$B^{y,t}_u=y$ for all $u\in[0,t]$. Otherwise, choose the
largest $T^y_j=s<t$, and put $B^{y,t}_u=y$ for $u\in[0,t-s)$
and $B^{y,t}_{t-s}=x+X_{x,j}$. Continue in this fashion to
complete the construction of $B^{y}_u,u\in[0,t]$. Note that
each $B^{y,t}_u$ is a rate one random walk with step
distribution $p(\cdot)$ and that the walks coalesce when
they meet: if $B^{x,t}_u= B^{y,t}_u$ for some $u\in[0,t]$,
then $B^{x,t}_s=B^{y,t}_s$ for all $u\le s\le t$. On the
event that no $*$-arrow is encountered along the path
$B^{x,t}_\cdot$, that is, $(z,T^{*, z}_n)\ne
(B^{x,t}_{t-u},t-u)$ for all $z,n$ and $0\le u\le t$, then
%
%
\begin{equation}
\label{vmdual} \xi^\vep_t(x) = \xi^\vep_0
\bigl(B^{x,t}_t \bigr) \qquad \forall\xi^\vep_0
\in\{0,1\}^{\mathbb{Z}^d}.
\end{equation}

%
\begin{lem}\label{vmirred} Fix $t>0$, distinct
$y_0,y_1\in{\mathbb{Z}^d}$ and finite disjoint
$B_0,B_1\subset{\mathbb{Z}^d}$. Then
there exists a finite
$\Lambda=\Lambda(y_0,y_1,B_0,B_1)\subset{\mathbb{Z}^d}$ and a
$\cG(\Lambda\times[0,t])$-measurable event
$G=G(t,y_0,y_1,B_0,B_1)$ such that $P(G)>0$ and on $G$:
\begin{longlist}[(iii)]
\item[(i)] $T^{*,z}_1>t$ for all $z\in\Lambda$;

\item[(ii)] $B^{x,t}_u\in\Lambda$ for all $x\in B_0\cup B_1$, $u\in[0,t]$;

\item[(iii)] $B^{x,t}_t=y_i$ for all $x\in B_i$, $i=0,1$.
\end{longlist}
If $\xi^\vep_0(y_i)=i$, $i=0,1$, then on the event $G$,
$\xi^\vep_t(x)=i$ for all $x\in B_i$, $i=0,1$.
\end{lem}
\begin{pf} We reason as in the proof
of Lemma~\ref{suffindecomp}, but now working with the dual of $\xi
^\vep
$, using the fact that the
$B^{y,t}_u$ are independent, irreducible random walks as long as they
do not meet. There are sets $B'_0,B'_1$
which are far apart, each with widely separated points,
such that a sequence of
walk steps can move the walks from $B_0$ to $B'_0$ and $B_1$
to $B'_1$ without collisions. If $B_0'$ and $B_1'$ are
sufficiently far apart, then by irreducibility there is a
sequence of steps resulting in the walks from $B_0'$
coalescing at some site $y'_0$, the walks from $B_1'$
coalescing at some site $y'_1$, all without collisions between the two
collections of walks, and with
$y'_0$ and $y'_1$ far apart.
Now by moving one walk at a time it is possible to prescribe a set of walk
steps which take the two walks from $y_0$ and $y_1$ to $y'_0$ and
$y'_1$, respectively,
without collisions between the two walks. By reversing these steps
(recall $p$ is symmetric)
we can therefore have the above walks follow steps which will take them
from $y'_0$ and $y'_1$ to
$y_0$ and $y_1$, respectively, without collisions. In this way we can
prescribe walk steps which occur with positive probability and ensure
that $B_t^{x,t}=y_i$ for all $x\in B_i$.
Let $\Lambda$ be a finite set large enough to
contain all the positions of the walks in this process, and
let $G$ be the event that $T^{*,x}_1>t$ for all
$x\in\Lambda$, and such that the $T^{x}_n$ and $X_{x,n}$, $x\in
\Lambda$,
allow for the above prescribed sequence of walk steps to occur by time
$t$. Then $G$ has the desired properties,
and on this event,
$\xi^\vep_t(x)=\xi^\vep_0(B^{x,t}_t)$ for all $x\in B_0\cup B_1$ by
\eqref{vmdual}. Now
the fact that (iii) holds on $G$, implies the final conclusion by the
choice of $y_i$.
\end{pf}

In addition to Lemma~\ref{vmirred} we will need the
simpler fact that for any fixed $t>0$ and $z\in{\mathbb{Z}^d}$,
%
%
\begin{equation}
\label{02z} \inf_{\xi^\vep_0\dvtx\xi^\vep_0(0)=1}P \bigl(\xi^\vep_t(z)=1
\bigr)>0.
\end{equation}
This is clear because there is a sequence of random walk steps
leading from $0$ to $z$, and there is positive probability
that the walk makes these steps before time $t$ and that no
other transitions occur at any site in the sequence.

%
\begin{rem} It is clear that the above holds equally well for
voter model perturbations in $d=2$.
\end{rem}

\section{A complete convergence theorem for cancellative
systems}\label{secflip}

To make effective use of annihilating duality we will need
to know that for large $t$, if $\xi_t\ne\zero,\one$ and
finite $A\subset{\mathbb{Z}^d}$ is large, then there will be many sites
in $\xi_t\cap A$ which can flip values in a fixed time
interval, and that the probability there will be an odd
number of these flips is close to $1/2$. For $x\in{\mathbb{Z}^d}$ and
$A\subset{\mathbb{Z}^d}$ define
\[
A(x,\xi) = \bigl\{y\in A\dvtx\xi(y)=1 \mbox{ and }\xi(y+x)=0 \bigr\}.
\]
The conditions we will use are: there exists $x_0\in{\mathbb{Z}^d}$
such that
%
%
\begin{equation}
\label{flip2}\qquad  \lim_{K\to\infty}\mathop{\sup_{A\subset{\mathbb
{Z}^d}}}_{|A|\ge K}
\limsup_{t\to
\infty} P \bigl(|\xi_{t}|>0 \mbox{ and }
A(x_0,\xi_t) = \varnothing\bigr)=0\qquad \mbox{if }|
\widehat\xi_0|=\infty
\end{equation}
and
%
%
\begin{equation}
\label{oddgoal} \lim_{K\to\infty}\mathop{\sup_{A\in Y,\xi_0\in\{
0,1\}^{\mathbb{Z}^d}\dvtx}}_{
|A(x_0,\xi
_0)|\ge K}
\bigl|P\bigl(|\xi_1\cap A|\mbox{ is odd}\bigr)-\tfrac12 \bigr| =0.
\end{equation}
We will verify in Lemmas~\ref{lemodd2} and \ref{lemflip2}
below that our voter model perturbations have these
properties for all sufficiently small $\vep$, but first we
will show how they are used along with \eqref{eqH}
to obtain complete convergence of $\xi_t$. Recall $\nu_{1/2}$ is the
translation invariant stationary measure in \eqref{nuhalf1}.

%
\begin{prop}\label{propCCT} Let $\xi_t$ be a translation
invariant cancellative spin-flip system with rate function
$c(x,\xi)$ satisfying \eqref{traps}--\eqref{q02},
\eqref{flip2}, and \eqref{oddgoal}. Let $\zeta_t$ be the
annihilating dual with $Q$-matrix given in
\eqref{qmatrix} and assume that \eqref{eqH} holds.
Then $\nu_{1/2}$
satisfies \eqref{coexist}, and if $|\widehat\xi_0|=\infty$ then
%
%
\begin{equation}
\label{eqgoal} \xi_{t} \To\beta_0(\xi_0)
\delta_0 + \bigl(1-\beta_0(\xi_0) \bigr)
\nu_{1/2} \qquad\mbox{as }t\to\infty.
\end{equation}
\end{prop}

\begin{pf}
We start with some preliminary facts. First,
\eqref{flip2} implies that for any
$m<\infty$ and $\xi_0\in\{0,1\}^{\mathbb{Z}^d}$ with
$|\widehat\xi_0|=\infty$,
%
%
\begin{equation}
\label{flip2a} \lim_{K\to\infty}\mathop{\sup_{A\subset{\mathbb
{Z}^d}}}_{|A|\ge K}
\limsup_{t\to
\infty} P \bigl(|\xi_{t}|>0 \mbox{ and }\bigl |
A(x_0,\xi_t)\bigr|<m \bigr)=0.
\end{equation}
This is because
$|A|\ge mK$ implies $A$ can be written as the disjoint union of
sets $A_1,\ldots,A_m$ with each $|A_i|\ge K$, and
\[
\bigl\{|\xi_{t}|>0\mbox{ and } \bigl| A(x_0,
\xi_t)\bigr|<m \bigr\} \subset\bigcup_{i=1}^m
\bigl\{|\xi_{t}|>0\mbox{ and } \bigl| A_i(x_0,
\xi_t)\bigr|=0 \bigr\}.
\]
Applying \eqref{flip2} we obtain \eqref{flip2a}.

Next, we need a slight upgrade of the basic duality equation. As
shown in \cite{G79}, \eqref{eqduality2} can be extended by
applying the Markov property of $\xi_t$ at a time $v<t$. If the
processes $\xi_t$ and $\zeta_t$ are independent, then for
all $u,v\ge0$,
%
%
\begin{equation}
\label{eqduality3} P\bigl(|\xi_{v+u} \cap\zeta_0| \mbox{ is
odd }\bigr)= P\bigl(|\xi_v \cap\zeta_u| \mbox{ is odd }\bigr).
\end{equation}

Let $\nu_{1/2}$ be defined by \eqref{nuhalf2}. Since we are
assuming
$|\hxi_0|=\infty$, we have $\beta_1(\xi_0)=0$ by \eqref{betanonzero}.
In view of \eqref{nuhalf2} and $\delta_{\zero}(|\xi\cap A|\mbox{
is odd})=0$, to prove \eqref{eqgoal} it suffices
to prove [recall \eqref{odddet}] that for fixed $A\in Y$,
%
%
\begin{equation}
\label{goalflip} \lim_{t\to\infty} P \bigl(|\xi_{t}\cap A|
\mbox{ is odd} \bigr) =\frac12 \beta_\infty(\xi_0)P \bigl(
\zeta^A_t\ne\varnothing\ \forall t\ge0 \bigr).
\end{equation}
Fix $\vep>0$. By \eqref{oddgoal} there exists
$K_1<\infty$ such that if $B\in Y$
and $|B(x_0,\xi_0)|\ge K_1$, then
%
%
\begin{equation}
\label{odd1} \bigl| P\bigl(|\xi_1\cap B|\mbox{ is odd}\bigr)-\tfrac12 \bigr| <\vep.
\end{equation}
By \eqref{flip2a}, there exists $K_2<\infty$ and $s_0<\infty$
such that if
$|B|\ge K_2$ and $s\ge s_0$, then
%
%
\begin{equation}
\label{flip3} P \bigl(\xi_s\ne\varnothing\mbox{ and } \bigl|
B(x_0,\xi_s)\bigr|<K_1 \bigr)< \vep.
\end{equation}
By \eqref{eqH} we can choose $T=T(A,K_2)<\infty$ large enough so
that
%
%
\begin{equation}
\label{dualbig1} P \bigl(0<\bigl|\zeta^A_T\bigr|\le
K_2 \bigr) < \vep.
\end{equation}

%
\begin{figure}

\includegraphics{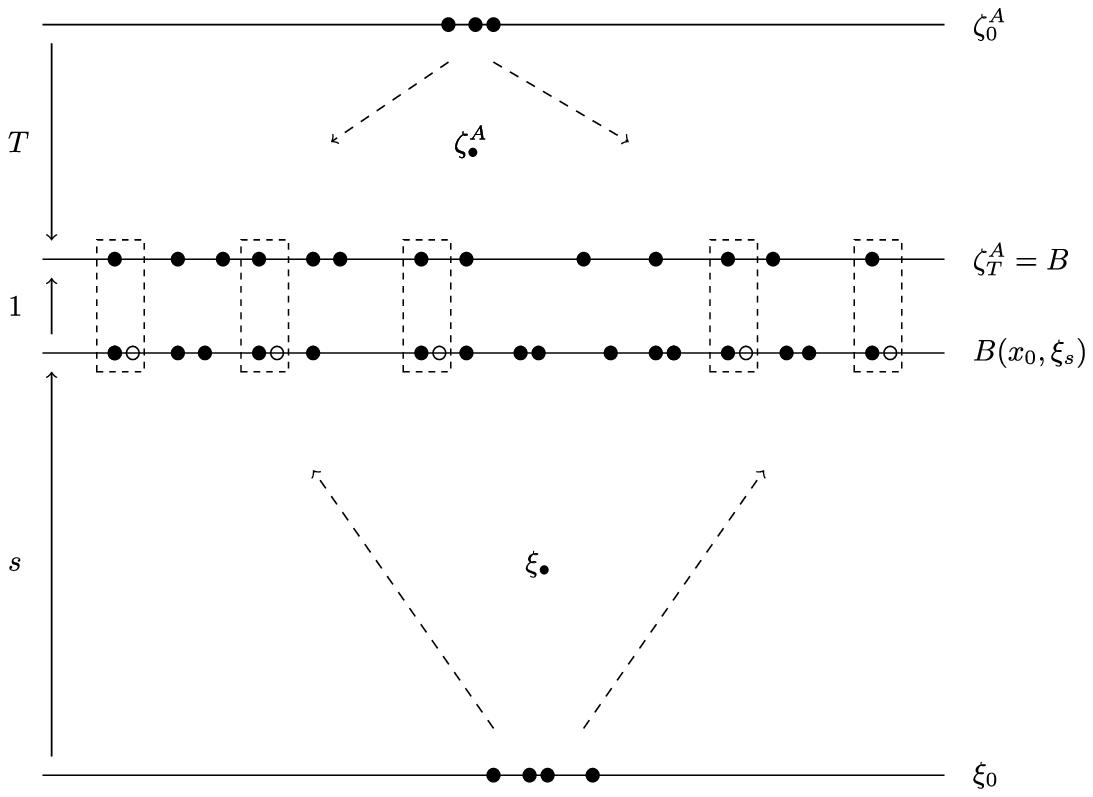}

\caption{$P(|\xi_{s+1}\cap\zeta^A_T|\mbox{ is odd})\approx
\frac12P(\xi_s\ne\varnothing)P(\zeta^A_T\ne\varnothing)$.}\label{fig1}
\end{figure}

For $t>1+T+s_0$ let $s=t-(1+T)$ and put $u=T$ and
$v=s+1$ in \eqref{eqduality3}. Then $P(|\xi_{t}\cap A|\mbox{ is
odd}) = P(|\xi_{s+1}\cap\zeta^A_T|\mbox{ is odd})$,
where $\xi_t$ and $\zeta^A_t$ are independent. (At this point the
reader may
want to consult Figure \ref{fig1} and Remark~\ref{Figdesc} below.) Making use of
the Markov property of $\xi_t$, we obtain
\begin{eqnarray*}
&&P\bigl(|\xi_{t}\cap A|\mbox{ is odd}\bigr)- \frac12 P(\xi_s\ne
\varnothing)P \bigl(\zeta^A_T\ne\varnothing\bigr)
\\
&&\qquad= \sum_{B\ne\varnothing}P \bigl(\zeta^A_T=B
\bigr) \biggl[ P\bigl(|\xi_{s+1}\cap B|\mbox{ is odd}\bigr)-\frac12 P(
\xi_s\ne\varnothing) \biggr]
\\
&&\qquad = \sum_{B\ne\varnothing}P \bigl(\zeta^A_T=B
\bigr) E \biggl[ \biggl(E_{\xi_s}\bigl(|\xi_1\cap B|\mbox{ is
odd}\bigr)-\frac12 \biggr) 1\{\xi_s\ne\varnothing\} \biggr].
\end{eqnarray*}
By \eqref{dualbig1},
%
%
\begin{eqnarray}
\label{decomp2}\qquad &&\biggl|P\bigl(|\xi_{t}\cap A|\mbox{ is odd}\bigr)- \frac12
P_\xi(\xi_s\ne\varnothing)P \bigl(\zeta^A_T
\ne\varnothing\bigr)\biggr|
\nonumber
\\[-8pt]
\\[-8pt]
\nonumber
&&\qquad< \vep+ \sum_{|B|> K_2}P \bigl(\zeta^A_T=B
\bigr) E \biggl[ \biggl|E_{\xi_s}\bigl(|\xi_1\cap B|\mbox{ is odd}\bigr)-
\frac12 \biggr| 1\{\xi_s\ne\varnothing\} \biggr].
\end{eqnarray}
By \eqref{flip3}, since $s>s_0$, each expectation in the last
sum is bounded above by
%
%
\begin{equation}
\label{decomp3} \vep+ E \bigl[ \bigl| E_{\xi_s}\bigl(|\xi_1\cap B|
\mbox{ is odd}\bigr)-\tfrac12 \bigr| 1 \bigl\{B(x_0,\xi_s)\ge
K_1 \bigr\} \bigr].
\end{equation}
Applying the bound \eqref{odd1} in this last expression,
and then combining \eqref{decomp2} and~\eqref{decomp3}
we obtain
\[
\bigl|P\bigl(|\xi_{t}\cap A|\mbox{ is odd}\bigr)- \tfrac12 P(\xi_s\ne
\varnothing)P \bigl(\zeta^A_T\ne\varnothing\bigr)\bigr| < 3
\vep.
\]
Let $t$ (and hence $s$) tend to infinity, and then $T$ tend to
infinity above to complete the proof of
\eqref{goalflip} and hence \eqref{eqgoal}.

Finally, let $|\xi_0|=|\hxi_0|=\infty$. Then
$\beta_0(\xi_0)=0$ by \eqref{betanonzero}, so \eqref{eqgoal} and
\eqref{flip2a}
imply that for any finite $m$, $\nu_{1/2}(|\xi|\ge m,
|\hxi|\ge m) =1$, and this implies coexistence.
\end{pf}

%
\begin{rem}\label{Figdesc} Figure~\ref{fig1} above gives a graphical view
of the
above argument. Time runs up for $\xi$ and down for
$\zeta$. Conditional on $\xi_s\ne\varnothing$ and
$\zeta^A_T\ne\varnothing$, \eqref{eqH} guarantees that
$B=\zeta^A_T$ is large, and \eqref{eqduality3} guarantees
$B(x_0,\xi_s)$ is large. In the dashed boxes in Figure \ref{fig1},
the $(\bullet\circ)$ pairs indicate the locations
$(x,x+x_0)$, $x\in B(x_0,\xi_s)$. Finally,
\eqref{oddgoal} now guarantees that $|\xi_{s+1}\cap
B|$ will be odd with probability approximately $\frac12$.
\end{rem}

Verification of \eqref{flip2} for our voter model
perturbations requires a comparison with oriented
percolation which we will save for the next section. Here we
present a proof of \eqref{oddgoal}, based heavily on ideas
from \cite{BDD}. See Lemma~7 of \cite{SS} for a purely cancellative version
of this result.

%
\begin{lem}\label{lemodd2}
If $\xi^\vep_\cdot$ is
a voter model perturbation,
then there exists $\vep_1>0$ and $x_0\in{\mathbb{Z}^d}$ such that
\eqref
{oddgoal} holds
for $\xi^\vep_\cdot$ if $\vep<\vep_1$.
\end{lem}

\begin{pf} Fix any $x_0$ with $p(x_0)>0$. We will prove
that if $\delta>0$, then there exists $K$ such that if
$|A(x_0,\xi_0)|\ge K$, then
%
%
\begin{equation}
\label{oddgoalep} \bigl| P \bigl(\bigl|\xi^\vep_1\cap A\bigr|\mbox{ is
odd} \bigr) -\tfrac12 \bigr| <\delta.
\end{equation}

Using the graphical construction of $\xi^\vep_t$ described
in Section~\ref{secirred}, we define a
version of the ``almost isolated sites'' of \cite{BDD}.
First we give the informal definition. For $x\in{\mathbb{Z}^d}$,
let $U(x)$ be the indicator of the
event that during the time period $[0,1]$, no change can
occur at site $x+x_0$ and no change can occur
at $x$ except possibly due to a (first) voter arrow directed from
$x$ to $x+x_0$. Let $V(x)$ be the indicator of
the event that during the time period $[0,1]$ no site
$y$ outside $\{x,x+x_0\}$ can change due to the value at $x$. More
formally, for $y\in{\mathbb{Z}^d}$ and $A\in Y$ with $|A|\le N_0$, define
\begin{eqnarray*}
\tau(y,A)& =& \min\bigl\{T^y_n\dvtx A=
\{X_{y,n} \}, n\in\NN\bigr\}
\\
&&{}\wedge\min\bigl\{T^{*,y}_n\dvtx A= \bigl
\{Z^1_{y,n},\ldots,Z^{N_0}_{y,n} \bigr\},
n \in\NN\bigr\},
\end{eqnarray*}
and $\tau(y)=\min\{\tau(y,A)\dvtx A\in Y$\}.
We can now define
\begin{eqnarray*}
U(x) &= &1 \bigl\{ \tau(x+x_0)>1, X_{x,1}=x_0,
T^x_2>1\mbox{ and }T^{*,x}_1>1 \bigr
\}\quad \mbox{ and}
\\
V(x) &= &1 \bigl\{ \tau(y,A)>1 \ \forall y\in{\mathbb{Z}^d}\setminus
\{x,x+x_0 \}\mbox{ and } A\in Y\colon x\in y+A \bigr\},
\end{eqnarray*}
and call $x$ almost isolated if $U(x)V(x)=1$.

By standard properties of Poisson processes,
%
%
\begin{equation}
\label{Poissonprop}\tau(x,A)\mbox{ and } \tau(y,B)\mbox{ are independent
whenever }x\ne y\mbox{ or }A\ne B.
\end{equation}
We also define
\[
\nu(A) = P \bigl(\{X_{0,1}\}=A \bigr) + P \bigl( \bigl
\{Z^1_{0,1},\ldots,Z^{N_0}_{0,1} \bigr\}=A
\bigr),
\]
and observe that $\nu(A)=0$ if $|A|>N_0$, and $\sum_{A\in
Y}\nu(A)= 2$. Use the fact that $\{T^0_n\dvtx\{X_{0,n}\}=A\}$ are the
points of a Poisson
point process with rate $P(\{X_{0,1}\}=A)$, and $\{T^{*,0}_n\dvtx\{
Z^1_{0,n},\ldots,Z_{0,n}^{N_0}\}=A\}$ are the points of an independent
Poisson point process with rate $\bar c\vep^2P(A=\{Z^1_{0,n},\ldots,Z_{0,n}^{N_0}\})$ to conclude that
%
%
\begin{eqnarray}
\label{rangeest}\qquad P \bigl(\tau(y,A)>1 \bigr)&=&\exp\bigl(-P \bigl(
\{X_{0,1} \}=A \bigr)-\bar c\vep^2P \bigl( \bigl\{
Z^1_{0,n}, \ldots,Z_{0,n}^{N_0} \bigr\}=A
\bigr) \bigr)
\nonumber
\\[-8pt]
\\[-8pt]
\nonumber
& \ge& e^{-(1+\bar c)\nu(A)}.
\end{eqnarray}

For
each $x\in{\mathbb{Z}^d}$, the variables $U(x),V(x)$ are independent
[this much is
clear from \eqref{Poissonprop}],
and we claim that $u_0=E(U(x))$ and $v_0=E(V(x))$ are positive
uniformly in $\vep$. To check this for $v_0$ we apply \eqref{rangeest}
to get
\begin{eqnarray*}
v_0 &\ge&\exp\biggl(-(1+\bar c)\sum_{A\in Y}
\sum_{y\in{\mathbb{Z}^d}\setminus
\{x\}} \nu(A) 1\{x-y\in A\} \biggr)
\\
& \ge&\exp\biggl(-(1+\bar c)\sum_{A\in Y}|A|\nu(A)
\biggr),
\end{eqnarray*}
which is positive, uniformly in $\vep\le\vep_0$.
For $u_0$, we have by the choice of $x_0$,
\[
u_0\ge\exp\biggl\{-(1+\bar c)\sum_{A\in Y}
\nu(A) \biggr\} p(x_0)P \bigl(T^0_2>1
\bigr)P \bigl(T_1^{*,0}>1 \bigr)>0.
\]
If $w_0=w_0(\vep)=u_0v_0$, then we have verified that
%
%
\begin{equation}
\label{gammalbnd}\gamma=\min\bigl\{w_0(\vep)\dvtx0<\vep\le\vep
_0 \bigr\}>0.
\end{equation}

Now suppose $\cY=\{y_1,\ldots,y_J\}\subset{\mathbb{Z}^d}$ and
$|y_i-y_j|>2|x_0|$ for
$i\ne j$. Then $U(y_1),\ldots,U(y_J)$ are independent, but
$V(y_1),\ldots,V(y_J)$ are not. Nevertheless, we claim
they are almost independent if all $|y_i-y_j|$, $i\ne
j$, are
large, and hence if we let $W(y_i)=U(y_i)V(y_i)$, then
$W(y_1),\ldots, W(y_J)$ are almost independent.
More precisely, we claim that for any $J\ge2$
and $a_i\in\{0,1\}$, $1\le i\le J$,
%
%
\begin{eqnarray}
\label{ai}\qquad && \lim_{n\to\infty} \mathop{\sup_{\cY=\{y_1,\ldots,y_J\},}}_{|y_i-y_j|\ge n
\ \forall i\ne j}
\Biggl|P \bigl( W(y_i)=a_i \ \forall1\le i\le J \bigr) -\prod
_{i=1}^JP \bigl(W(y_i)=a_i
\bigr) \Biggr|
\nonumber
\\[-8pt]
\\[-8pt]
\nonumber
&&\qquad = 0.
\end{eqnarray}
For the time being, let us suppose this fact.

Given $J$ and $\cY=\{y_1,\ldots,y_J\}$, let
$S(J,\cY)=\sum_{y\in\cY}W(y)$. Then \eqref{ai} implies
$S(J,\cY)$ is approximately binomial if the $y_i$ are well
separated. That is, if $\cB(J,w_0)$ is a binomial random
variable with parameters $J$ and $w_0$, and
\[
\Delta(J,\cY,k)= \bigl|P \bigl(S(J,\cY)=k \bigr) - P \bigl(\cB(J,w_0)=k
\bigr) \bigr|,
\]
then \eqref{ai} implies that for $k=0,\ldots,J$,
%
%
\begin{equation}
\label{spreadout} \lim_{n\to\infty}\mathop{\sup_{\cY=(y_1,\ldots,y_J)}}_{
|y_i-y_j|\ge n, i\ne j}
\Delta(J,\cY,k) =0.
\end{equation}

Now fix $\delta>0$. A short calculation
shows that
%
%
\begin{equation}
\label{mireq} p_0=P \bigl(T^x_1>1|U(x)=1
\bigr)=P \bigl(T_1^x>1|T^x_2>1
\bigr)=\tfrac{1}{2}.
\end{equation}
By \eqref{gammalbnd} and
\eqref{spreadout},
we may choose $J=J(\delta)$ such that
%
%
\begin{equation}
\label{oddgamma} (1-\gamma)^J <\delta,
\end{equation}
and then $n=n(J,\delta)$ so that
for all $\cY=\{y_1,\ldots,y_J\}$ with $|y_i-y_j|\ge n$ for
$i\ne j$,
%
%
\begin{equation}
\label{binomialapprox} \Delta(J,\cY,0) <\delta.
\end{equation}
Given $J$ and $n$, it is easy to see that there exists $K=K(J,n)$ such
that if $B\subset{\mathbb{Z}^d}$ and $|B|\ge K$, then $B$ must contain
some $\cY=\{y_1,\ldots,y_J\}$ such that $|y_i-y_j|\ge n$ for $i\ne j$.

Now suppose that $|A(x_0,\xi^\vep_0)|\ge K$ and
$\cY=\{y_1,\ldots, y_J\}\subset A(x_0,\xi^\vep_0)$ with
$|y_i-y_j|\ge n$
for all $i\neq j$. Let
$\cI$ be the set of $y_j$ with $W(y_j)=1$, so that
$|\cI|=S(J,\cY)$.
Let $\cG$ be the $\sigma$-field generated by
%
%
\begin{eqnarray}
\label{list} &&\bigl\{1(y_j\in\cI)\dvtx j=1,\ldots,J \bigr\}
\nonumber
\\[-8pt]
\\[-8pt]
\nonumber
&&\qquad{}\cup\bigl\{1 \bigl\{x\in\cI^c \bigr\} \bigl
(T_n^x,T^{*,x}_n,X_{x,n},Z_{x,n},U_{x,n}
\dvtx x\in{\mathbb{Z}^d},n\ge1 \bigr) \bigr\}.
\end{eqnarray}
If $g_j=1\{T^{y_j}_1>1\}$, then conditional on $\cG$,
%
%
\begin{equation}
\label{Bin2} \{g_j\dvtx y_j\in\cI\}\mbox{ are i.i.d.
Bernoulli rv's with mean $p_0=\tfrac{1}{2}$},
\end{equation}
and $X=\sum_j 1\{y_j\in\cI\}g_j$ is binomial with parameters $(|\cI
|,p_0=\frac{1}{2})$.
This is easily checked by conditioning on the $\cG$-measurable set
$\cI
$ and using \eqref{mireq}.

Let $h=\sum_{x\in A}1\{x\notin\cI)\xi^\vep_1(x)$. Then at time 1 we
have the
decomposition
%
%
\begin{equation}
\label{decomp1} \bigl|\xi^\vep_1\cap A\bigr|=h + X,
\end{equation}
where we have used the fact that for $y_j\in\cI$, $\xi^\vep_s(y_j)$
will flip from a $1=\xi^\vep_0(y_j)$ to a $0=\xi^\vep_0(y_j+x_0)$
during the time interval
$[0,1]$ if and only if $g_j=0$. Since $h$ is $\cG$-measurable,
\begin{eqnarray}
P \bigl( \bigl|\xi^\vep_1\cap A\bigr| \mbox{ is odd}\mid\cG\bigr)
( \omega) = P \bigl( X =1-h(\omega) \bmod2 \mid\cG\bigr) (\omega)=
\tfrac{1}{2} \nonumber\\
\eqntext{\mbox{a.s. on }\bigl\{|\cI|>0\bigr\},}
\end{eqnarray}
the last by an elementary binomial calculation and the fact that
conditional on $\cG$,
$X$ is binomial with parameters $(|\cI|,\frac{1}{2})$.
Take
expectations in the above and use~\eqref{binomialapprox} and then
\eqref{oddgamma} to conclude that
%
%
\begin{eqnarray}
\label{odd2} \bigl|P \bigl(\bigl| \xi^\vep_1\cap A\bigr| \mbox{ is odd}
\bigr) - \tfrac12\bigr| &\le& P\bigl(|\cI|=0\bigr)
\nonumber\\
&\le&\delta+P \bigl(\cB(J,\gamma)=0 \bigr)
\\
\nonumber
&\le&\delta+(1-\gamma)^J\le2\delta.
\end{eqnarray}
This yields inequality \eqref{oddgoalep}.

It remains to verify \eqref{ai}. This is easy if $p$ and $q_Z$ have
finite support; see Remark~\ref{finsuppeasy} below. In general, the
idea is to write
%
%
\begin{equation}
\label{Vprod} V(y_i)=V_1(y_i)V_2(y_i)
V_3(y_i),
\end{equation}
where $V_1(y_1),\ldots, V_1(y_J)$ are independent and
independent of $U(y_1),\ldots,U(y_J)$, and
$V_2(y_1),V_3(y_1),\ldots,V_2(y_J),V_3(y_J)$ are all one with high
probability if the $y_i$ are sufficiently spread out.
Let $n\ge2|x_0|$ and $\cY=\{y_1,\ldots,y_J\}$ be given with
$|y_i-y_j|\ge n$, and let
$\cY_0=\{y_1+x_0,\ldots,y_J+x_0\}$. Note that
$y_1,y_1+x_0,\ldots,y_J,y_J+x_0$ are distinct. Define
\begin{eqnarray*}
V_1(y_i)&=& 1 \bigl\{ \tau(z,A)>1 \ \forall z\notin\cY
\cup\cY_0, A\in Y\dvtx(z+A)\cap\cY=\{y_i\} \bigr\},
\\
V_2(y_i)&=&1 \bigl\{\tau(z,A)>1 \ \forall z\notin\cY\cup
\cY_0, A\in Y\dvtx z+A\supset\{y_i,y_j\}\mbox{ for some }j\ne i \bigr\},
\\
V_3(y_i)&= &1 \bigl\{\tau(z,A)>1 \ \forall z\in(\cY\cup
\cY_0)\setminus\{y_i,y_i+x_0
\}, A\in Y\dvtx y_i\in z+A \bigr\}.
\end{eqnarray*}
A bit of elementary logic shows that \eqref{Vprod} holds.
If a pair $(z,A)$ occurs in the definition of some $V_1(y_i)$, then it
cannot occur in any $V_1(y_j), j\ne i$, and hence
$V_1(y_1),\ldots,V_1(y_J)$ are independent, and also
independent of $U(y_1),\ldots,U(y_J)$. Therefore, to prove
\eqref{ai} it suffices to prove that
%
%
\begin{equation}
\label{ai2} \lim_{n\to\infty} \mathop{\sup_{\cY=\{y_1,\ldots,y_J\},}}_{|y_i-y_j|\ge n
\ \forall i\ne j}P
\bigl(V_2(y_i)V_3(y_i)\ne1
\bigr) = 0.
\end{equation}

We treat $V_3(y_i)$ first.
By \eqref{rangeest},
\begin{eqnarray*}
P \bigl(V_3(y_i)=1 \bigr) &\ge&\exp\biggl(-(1+\bar c)
\sum_{z\in(\cY\cup\cY_0)\setminus\{
y_i,y_i+x_0\}} \sum_{A\in Y}
\nu(A)1\{y_i\in z+A\} \biggr)
\\
&\ge&\exp\biggl(-2J(1+\bar c) \sum_{A\in Y}\nu(A)1
\bigl\{\operatorname{diam}(A)>n \bigr\} \biggr)
\\
&\to& 1\qquad \mbox{as }n\to\infty.
\end{eqnarray*}
To treat $V_2(y_i)$ we note that
if a pair $(z,A)$ occurs in the definition of $V_2(y)$, then
$\operatorname{diam}(A)\ge n$, so
\begin{eqnarray*}
\hspace*{-4pt}&&P \bigl(V_2(y_i)=1 \bigr)
\\
\hspace*{-4pt}&&\quad\ge\exp\biggl(-(1+\bar c)\sum_{B\subset\cY}\sum
_{A\in Y}\sum_{z\ne y_i} 1 \bigl
\{y_i\in B= (z+A)\cap\cY,\operatorname{diam}(A)\ge n \bigr\}\nu(A) \biggr).
\end{eqnarray*}
In the sum above, given $B\ni y_i$ there at most $|A|$ choices
for $z$ such that $(z+A)\cap\cY=B$. In fact, there are at most $|A|$ choices
of $z$ such that $y_i\in z+A$ as this implies $z\in y_i-A$.
Thus
\begin{eqnarray*}
P \bigl(V_2(y_i)=1 \bigr) &\ge&\exp\biggl(-(1+\bar c)
\sum_{B\subset\cY}\sum_{A\in Y} 1
\bigl\{\operatorname{diam}(A)\ge n \bigr\}|A|\nu(A) \biggr)
\\
&\ge&\exp\biggl(-(1+\bar c) 2^J\sum_{A\in Y}
1 \bigl\{\operatorname{diam}(A)\ge n \bigr\}|A|\nu(A) \biggr)
\\
&\to&1 \qquad\mbox{as }n\to\infty.
\end{eqnarray*}

This proves \eqref{ai2} and hence \eqref{ai}.
\end{pf}

%
\begin{rem}\label{finsuppeasy}Note that if $p(\cdot)$ and
$q_Z(\cdot
)$ have finite support, then the proof
simplifies somewhat because for large enough $n$ the left-hand side of
\eqref{ai} is zero. This is because the $A$'s arising in the definition
of $V(x)$ will have uniformly bounded diameter which will show that for
$|y_i-y_j|$ large, $V(y_1),\ldots,V(y_J)$ will be independent.
\end{rem}

\section{\texorpdfstring{Proof of Theorem~\protect\ref{thmCCTpert}}
{Proof of Theorem 1.2}}
\label{secthmproof}
To
prove Theorem~\ref{thmCCTpert} it will suffice, in view of
Proposition~\ref{propCCT}, Lemmas~\ref{lemH} and
\ref{lemodd2} and Remark~\ref{remHaltcond}, to prove that for small
enough~$\vep$,
both conditions \eqref{eqH3} and \eqref{flip2} hold
for $\xi^\vep_t$. The proof of \eqref{flip2} is given in
Lemma~\ref{lemflip2} below after first developing the
necessary oriented percolation machinery.
With \eqref{flip2} in hand the proof of \eqref{eqH3} is
then straightforward.

We suppose now that $\xi^\vep$ is a voter model perturbation
with rate function $c_\vep(x,\xi)$ and that all the
assumptions of Theorem~\ref{thmCCTpert} are in force.
We also assume that $\xi^{\vep}$ is constructed
using the Poisson processes $T^x_n,T^{*,x}_n$ and the
variables $X_{x,n},Z^i_{x,n},U_{x,n}$ as in
Section~\ref{secirred}. We assume
$|{\hat\xi}^\vep_0|=\infty$, so that by \eqref{betanonzero}
$\beta_1(\xi^\vep_0)=0$.
By the results of \cite{CDP11} for small $\vep$ we expect that when
$\xi
^\vep_t$ survives, there will be
blocks in space--time, in the graphical construction,
containing both 0's and 1's, which dominate a
super-critical oriented percolation. The percolation process
necessarily spreads out. So if $A\subset{\mathbb{Z}^d}$ is large,
eventually there will be many blocks
containing 0's and 1's near the sites of $A$ at times just
before $t$,
allowing for many independent tries to force $|A(\xi_t,x_0)|\ge1$.

Let $\Zde$ be the
set of $x\in{\mathbb{Z}^d}$ such that $\sum_iz_i$ is even. Let
$\cL=\{(x,n)\subset{\mathbb{Z}^d}\times\Z^+\dvtx\sum_i x_i
+n\mbox{ is
even}\}$. We equip $\cL$ with edges from $(x,n)$ to
$(x+e,n+1)$ and $(x-e,n+1)$ for all
$e\in\{e_1,\ldots,e_d\}$,
where $e_i$ is the $i$th unit basis vector. Given a family of Bernoulli random
variables $\theta(x,n), (x,n)\in\cL$, we define open paths
in $\cL$ using the $\theta(y,n)$ and the edges in $\cL$ in the
usual way. That is, a sequence of points $z_0,\ldots,z_n$ in $\cL$ is an
open path from $z_0$ to $z_n$ if and only if there is an edge from
$z_i$ to $z_{i+1}$ and $\theta(z_i)=1$ (in which case we say site $z_i$
is open) for $i=0,\ldots,n-1$. We will write
$(x,n)\to(y,m)$ to indicate there is an open path in $\cL$
from $(x,n)$ to $(y,m)$. Define the open cluster starting at
$(x,n)\in\cL$,
\[
\cC(x,n) = \bigl\{ (y,m)\in\cL\dvtx m\ge n \mbox{ and } (x,n)\to(y,m)
\mbox{ in
} \cL\bigr\}.
\]
For $(x,n)\in\cL$ let
$W^{(x,n)}_{m}=\{y\dvtx(x,n)\to(y,m)\}$, $m\ge n$.
We will write $W^0_n$ for
$W^{(0,0)}_n$. For
$k=1,\ldots,d$, say that $(x,n)\to_k
(y,m)$ if there is an open path from $(x,n)$ to $(y,m)$
using only edges of the form $(x,n)\to(x+e_k,n+1)$ or
$(x,n)\to(x-e_k,n+1)$. We define the
corresponding ``slab'' clusters $\cC_k(x,n)$ and processes
$W^{(x,n)}_{k,m}$ using these paths.
Clearly $\cC_k(x,n)\subset\cC(x,n)$ and $W^{(x,n)}_{k,m}\subset W^{(x,n)}_m$.
If $W_0\subset\Z^d$, let $W_m=\bigcup_{x\in W_0}W_m^{(x,0)}$.

%
\begin{lem}\label{perclem} Suppose the $\{\theta(z,n)\}$ are i.i.d., and
$1-\gamma=P(\theta(x,n)=1)\ge1-6^{-4}$. Then
%
%
\begin{equation}
\label{percolate} \rho_\infty=P \bigl(\bigl|\cC_1(0,0)\bigr|=\infty
\bigr)>0
\end{equation}
and
%
%
\begin{equation}
\label{0occupied} \lim_{K\to\infty} \mathop{\sup
_{A\subset
2{\mathbb{Z}^d}}}_{|A|\ge K} \limsup_{n\to\infty}P
\bigl(W^0_{2n}\ne\varnothing\mbox{ and }
W^0_{2n}\cap A=\varnothing\bigr) =0.
\end{equation}
\end{lem}
\begin{pf}
For \eqref{percolate}, see Theorem~A.1 (with $M=0$) in
\cite{Dur95}. The limit \eqref{0occupied} is known for
$d=1$, while the $d>1$ case is an immediate consequence of
the ``shape theorem'' for $W^0_{2n}$, the discrete time
analogue of the shape theorem for the contact process in
\cite{DG82}. Since this discrete time result does not
appear in the literature, we will give a direct proof of
\eqref{0occupied}, but for the sake of simplicity will
restrict ourselves to the $d=2$ case. We need the
following $d=1$ results, which we state using our ``slab''
notation,
%
%
\begin{equation}
\label{percdens} \exists\rho_1>0 \mbox{ such that}\qquad\liminf
_{n\to\infty} P \bigl((x,0)\in W^{0}_{1,2n} \bigr)
\ge\rho_1\qquad \mbox{for all } x\in2\Z,\hspace*{-35pt}
\end{equation}
and
for fixed $K_0\in\NN$,
%
%
\begin{equation}
\label{K0occupied1d1} \lim_{K\to\infty}\sup_{A\subset2\Z\times\{
0\}, |A|\ge K}
\limsup_{n\to
\infty} P \bigl(W^0_{1,2n}\ne
\varnothing,\bigl|W^0_{1,2n}\cap A \bigr|<K_0\bigr) = 0.
\end{equation}
These facts are easily derived using the methods in
\cite{Dur84}; see also Lemma~3.5 in~\cite{DN}, the Appendix
in \cite{DN91} and Section~2 of \cite{BN94}.

The idea of the proof of the $d=2$ case of
\eqref{0occupied} is the following. If $n$ is large, then on the event
$W^0_{2n}\ne\varnothing$ we can find, with high probability, a point
$z\in
W^0_{2k}$ for some small $k$ such that
$W^{(z,2k)}_{1,2m}\ne\varnothing$ for some large $m<n$. With
high probability $W^{(z,2k)}_{1,2m}$ will contain many
points $z'$ from which we can start independent ``$e_2$'' slab processes
$W^{(z',2m)}_{2,2n}$. Many of these will be large, providing
many independent chances for $W^{(z',2m)}_{2,2n}\cap A\ne\varnothing$,
forcing $W^0_{2n} \cap A\ne\varnothing$.

Here are the details. We may assume without loss of
generality that all sets $A$ considered here are finite. Fix
$\delta>0$, and choose positive integers $J_0,K_0$
satisfying $(1-\rho_\infty)^{J_0}<\delta$ and
$(1-\rho_1)^{K_0}<\delta$. By \eqref{K0occupied1d1} we can
choose a positive integer $K_1$ such that for all
$A\subset2\Z\times\{0\}$, $|A|\ge K_1$,
%
%
\begin{equation}
\label{K0occupied1d2} \limsup_{n\to\infty}P \bigl(W^0_{1,2n}
\ne\varnothing, \bigl|W^0_{1,2n}\cap A\bigr|<K_0 \bigr)
< \delta.
\end{equation}
For $x=(x_1,x_2)\in\Z^2$ and $A\subset\Z^2$ let $\pi_1x=(x_1,0)$,
$\pi_2x=(0,x_2)$ and $\pi_iA=\{\pi_ia\dvtx a\in A\}$, $i=1,2$.
Observe that at least one of the $|\pi_iA|\ge\sqrt{|A|}$.
We now fix any $A\subset2\Z^2$ with $|A|\ge K_1^2$, and suppose
$|\pi_1A|\ge K_1$. For convenience later in the argument,
fix any $A'\subset A$ such that $\pi_1A'=\pi_1A$ and
$\pi_1$ is one-to-one on $A'$.
By \eqref{K0occupied1d2} we may choose a positive integer
$n_1=n_1(A)$ such that
%
%
\begin{equation}
\label{K0occupied1d3} P \bigl( W^0_{1,2n}\ne\varnothing,
\bigl|W^0_{1,2n}\cap\pi_1 A'\bigr|<K_0
\bigr) < \delta\qquad\mbox{for all }n\ge n_1.
\end{equation}
We may increase $n_1$ if necessary so that
$P(|\cC_1(0,0)|<\infty, W^0_{1,2n_1}\ne\varnothing)<\delta$,
which implies that
%
%
\begin{equation}
\label{dieaftern1} P \bigl(W^0_{1,2n_1}\ne\varnothing,
W^0_{1,2n}=\varnothing\bigr) < \delta\qquad\mbox{for all }n\ge
n_1.
\end{equation}

Let $m(j)=2(j-1)n_1, j=1,2,\ldots,$ and define a random
sequence of points $z_1,z_2,\ldots$ as
follows. If $W^{0}_{m(j)}\ne\varnothing$, let $z_j$ be the
point in $W^{0}_{m(j)}$ closest to the origin, with some
convention in the case of ties. If
$W^{0}_{m(j)}=\varnothing$, put $z_j=0$. Define
\[
N = \inf\bigl\{j \dvtx z_j\in W^{0}_{m(j)}
\mbox{ and } W^{(z_j,m(j))}_{1,m(j+1)}\ne\varnothing\bigr\}.
\]
Since $P(W^0_{1,2n_1}=\varnothing)\le1- \rho_\infty$, the
Markov property implies
\begin{eqnarray*}
P \bigl(W^0_{m(j)}\ne\varnothing, N>j \bigr) &=& P
\bigl(W^0_{m(j)}\ne\varnothing, N>j-1, W^{(z_j,m(j))}_{1,m(j+1)}=
\varnothing\bigr)
\\
&\le&(1- \rho_{\infty}) P \bigl(W^0_{m(j)}\ne
\varnothing, N>j-1 \bigr).
\end{eqnarray*}
The above is at most $(1- \rho_{\infty}) P(W^0_{m(j-1)}\ne
\varnothing,
N>j-1)$, so iterating this, we get
%
%
\begin{equation}
\label{geombnd}P \bigl(W^0_{m(j)}\ne\varnothing, N>j
\bigr) \le(1-\rho_\infty)^j,
\end{equation}
and if $n>J_0n_1$, then
%
%
\begin{equation}
\label{geombound}\quad  P \bigl(W^0_{2n}\ne\varnothing,
N>J_0 \bigr) \le P \bigl(W^0_{m(J_0)}\ne
\varnothing, N>J_0 \bigr) \le(1-\rho_{\infty})^{J_0}<
\delta.
\end{equation}

We need a final preparatory inequality. Using
\eqref{K0occupied1d3} and the Markov property, for $n>n_1$
we have
\begin{eqnarray*}
&&P \bigl(W^0_{1,2n}\ne\varnothing, W^0_{2n}
\cap A=\varnothing\bigr)
\\
&&\qquad\le\delta+\sum_{B\subset A', |B|\ge K_0}P \bigl(W^0_{1,2n_1}
\cap\pi_1A'=\pi_1B \bigr) P \bigl(x\notin
W^{(\pi_1x,2n_1)}_{2,2n}\ \forall x\in B \bigr)
\\
&&\qquad\le\delta+\sum_{B\subset A', |B|\ge K_0}P \bigl(W^0_{1,2n_1}
\cap\pi_1A'=\pi_1B \bigr) \prod
_{x\in B}P \bigl(x\notin W^{(\pi_1x,2n_1)}_{2,2n} \bigr),
\end{eqnarray*}
the last step by independence of the slab processes. Thus,
employing \eqref{percdens},
%
%
\begin{equation}
\label{prelim} \limsup_{n\to\infty}P \bigl(W^0_{1,2n}
\ne\varnothing, W^0_{2n}\cap A=\varnothing\bigr)
\\
\le\delta+(1-\rho_1)^{K_0} <2\delta.
\end{equation}

We are ready for the final steps. For each $j\le J_0$ and
$n\ge J_0n_1$, by
the Markov property and \eqref{dieaftern1},
\begin{eqnarray*}
&&P \bigl(W^0_{2n}\ne\varnothing,W^0_{2n}
\cap A=\varnothing,N=j \bigr)
\\
&&\qquad\le P \bigl(W^0_{m(j)}\ne\varnothing,N>j-1,
W^{(z_j,m(j))}_{1,m(j+1)}\ne\varnothing, W^{(z_j,m(j))}_{2n}
\cap A=\varnothing\bigr)
\\
&&\qquad=\sum_{z}P \bigl(W^0_{m(j)}
\ne\varnothing, N>j-1,z_j=z \bigr)
\\
&&\qquad\qquad{} \times P \bigl(W^{(z,m(j))}_{1,m(j+1)}\ne\varnothing, W^{(z,m(j))}_{2n}
\cap A=\varnothing\bigr)
\\
&&\qquad\le\sum_{z}P \bigl(W^0_{m(j)}
\ne\varnothing, N>j-1,z_j=z \bigr)
\\
&&\hspace*{14pt}\qquad\quad{} \times\bigl(\delta+ P \bigl(W^{(z,m(j))}_{1,2n}\ne\varnothing,
W^{(z,m(j))}_{2n}\cap A=\varnothing\bigr) \bigr).
\end{eqnarray*}
Applying \eqref{prelim} and then \eqref{geombnd}, we obtain
\begin{eqnarray*}
\limsup_{n\to\infty} P \bigl(W^0_{2n}\ne
\varnothing,W^0_{2n}\cap A=\varnothing,N=j \bigr) &\le&3
\delta P \bigl(W^0_{m(j)}\ne\varnothing,N>j-1 \bigr)
\\
&\le&3\delta(1-\rho_\infty)^{j-1}.
\end{eqnarray*}
It follows that
\begin{eqnarray*}
&&\limsup_{n\to\infty} P \bigl(W^0_{2n}\ne
\varnothing,W^0_{2n}\cap A=\varnothing\bigr)
\\
&&\qquad \le\limsup_{n\to\infty} P \bigl(W^0_{2n}
\ne\varnothing,N>J_0 \bigr)+ 3\delta\sum
_{j=1}^{J_0}(1- \rho_\infty)^{j-1}
\\
&&\qquad\le\delta+ 3\delta/\rho_\infty
\end{eqnarray*}
by using \eqref{geombound} and summing the series.
This completes the proof.
\end{pf}
Now we follow \cite{Dur95} and Section~6 of \cite{CDP11} in describing
a setup which
connects our spin-flip systems with the percolation process
defined above. Let $K,L,T$ be finite positive constants with $K, L\in
\NN
$, let
$r=\frac{1}{16d}$, $Q_\vep=[0,\lceil\vep^{r-1}\rceil]^d\cap\Z^d$ and
$Q(L)=[-L,L]^d$.
We define a set $H$ of configurations in $\{0,1\}^{\Z^d}$ to be an
unscaled version of the set of
configurations in $\{0,1\}^{\vep\Z^d}$ of the same name in Section~6 of
\cite{CDP11}, that is,
\begin{eqnarray*}
&&H= \biggl\{\xi\in\{0,1\}^{\Z^d}\dvtx|Q_\vep|^{-1}
\sum_{y\in Q_\vep}\xi(x+y)\in I^*\\
&&\hspace*{33pt}\mbox{for all }x\in Q(L)
\cap\bigl( \bigl[0, \bigl\lceil\vep^{r-1} \bigr\rceil
\bigr]^d \cap\Z^d \bigr) \biggr\}.
\end{eqnarray*}
Here $I^*$ is a particular closed subinterval of $(0,1)$; it is
$I^*_\eta$ in the notation of Section~6 in \cite{CDP11}. The key
property we will need of $H$ is
%
%
\begin{equation}
\label{Hprop} \mbox{for each $\xi\in H$ there are $y_0,y_1
\in Q(L)\cap\Z^d$ s.t. $\xi(y_i)=i$ for $i=0,1$.}\hspace*{-35pt}
\end{equation}
This is immediate from the definition and the fact that $I^*$ is a
closed subinterval of $(0,1)$.
For $z\in\Z^d$, let $\sigma_z\dvtx\{0,1\}^{\Z^d}\to\{0,1\}^{\Z
^d}$ be the
translation map, $\sigma_z(\xi)(x)=\xi(x+z)$, and let $0<\gamma'<1$.
Recall from Section~\ref{secirred} that for $R\subset\R^d$, $\cG
(R\times[0,T])$ is the $\sigma$-field generated by the points of the
graphical construction in the space--time region $R\times[0,T]$. For
each $\xi\in H$,
$G_\xi$ will denote an event such that:
\begin{longlist}[(iii)]
\item[(i)]$G_\xi$ is $\cG([-KL,KL]^d\times[0,T])$-measurable;
\item[(ii)] if $\xi_0=\xi\in H$, then on $G_\xi$,
$\xi_T\in\sigma_{Le}H$ for all
$e\in\{e_1,-e_1,\ldots,e_d,-e_{d}\}$;
\item[(iii)]$P(G_\xi)\ge1-\gamma'$ for all
$\xi\in H$.
\end{longlist}

We are now in a position to quote the facts we need from
Section~6 of \cite{CDP11}, which depend heavily on our
assumption $f'(0)>0$ [and by symmetry $f'(1)=f'(0)>0$]. This allows
us to use Proposition~1.6 of \cite{CDP11} to show that
Assumption 1 of that reference is in force and so by a
minor modification of Lemma~6.3 of \cite{CDP11} we have the following.

%
\begin{lem}\label{vmpertspercolate}
For any $\gamma'\in(0,1)$ there exists
$\vep_1>0$ and finite $K\in\NN$ such that for all
$0<\vep<\vep_1$ there exist $L,T,\{G_\xi,\xi\in
H\}$, all depending on $\vep$, satisfying the basic setup
given above.
\end{lem}

Lemma 6.3 of \cite{CDP11} deals with a rescaled process on
the scaled lattice $\vep\Z^d$ but here we have absorbed the
scaling parameters into our constants $T$ and $L$ and then
shifted $L$ slightly so that it is a natural number. In fact
$L$ will be of the form $\lceil c
\vep^{-1}\log(1/\vep)\rceil$.

Given $\xi=\xi_0\in\{0,1\}^{\mathbb{Z}^d}$
we define
%
%
\begin{equation}
\label{Vndef} V_n= \bigl\{x\dvtx(x,n)\in\cL\mbox{ and }
\sigma_{-Lx} \xi_{nT}\in H \bigr\}.
\end{equation}
Note that
$V_n=\varnothing$ and $V_{n+1}\ne\varnothing$ is possible.
Theorem~A.4 of \cite{Dur95} and its proof imply that
there are $\{0,1\}$-valued random variables $\{\theta'(z,n)\dvtx
(z,n)\in
\cL
\}$ so that
if $\{{W_m'}^{(x,n)}\dvtx m\in\Z_+, (x,n)\in\cL\}$ and $\{\cC
'(z,n)\dvtx(z,n)\in
\cL\}$ are constructed from $\{\theta'(z,n)\}$ as above, then
%
%
\begin{equation}
\label{VWprime} \mbox{if $x\in V_n$, then ${W'_m}^{(x,n)}
\subset V_m$ for all $m\ge n$,}
\end{equation}
and $\{W'_n\}$ is a $2K$-dependent oriented percolation process, that is,
%
%
\begin{eqnarray}
\label{modMdep} P \bigl(\theta'(z_k,n_k)
= 1 \mid\theta'(z_j,n_j),j< k \bigr)
\ge1- \gamma'
\end{eqnarray}
whenever $(z_j,n_j), 1\le j\le k$,
satisfy $n_j<n_k$, or $n_j=n_k$ and $|z'_j-z'_k|> 2K$, for all $j<k$.
The Markov property of $\xi^\vep$ allows us to only require $n_j<n_k$
as opposed to $n_k-n_j>2K$
in the above, as in Section~6 of \cite{CDP11}.

Let $\Delta=(2K+1)^{d+1}$. By Theorem~B26 of \cite{Lig99},
modified as in Lemma~5.1 of~\cite{CDP11},
if $\gamma'$ (in Lemma~\ref{percolate}) is taken small enough so
that
$1-\gamma= (1-(\gamma')^{1/\Delta})^2 \ge1/4$,
then the $\theta'(z,n)$ can be coupled with
i.i.d. Bernoulli variables $\theta(z,n)$ such that
%
%
\begin{eqnarray}
\label{thetathetaprime} %
\theta(z,n) &\le&\theta'(z,n)
\qquad\mbox{for all }(z,n)\in\cL\quad \mbox{and}
\nonumber
\\[-8pt]
\\[-8pt]
\nonumber
P \bigl(\theta(z,n)=1 \bigr)&=&1-\gamma.
\end{eqnarray}
(The simpler condition on $\gamma$ and $\gamma'$ in Theorem B26 of
\cite{Lig99}
and above in fact follows from that in \cite{LSS} and Lemma~5.1 of
\cite
{CDP11} by
some arithmetic, and the explicit value of $\Delta$ comes from the fact
that we are now
working on $\Z^d$.)
If the coupling part of \eqref{thetathetaprime}
holds, then $W_n\subset V_n$ for all $n$,
and \eqref{VWprime} implies
%
%
\begin{equation}
\label{percext} x\in V_n \mbox{ implies } W^{(x,n)}_m
\subset V_m \mbox{ for all } m\ge n.
\end{equation}
Now choose $\gamma'$ small enough in Lemma~\ref{vmpertspercolate} so that
%
%
\begin{equation}
\label{LSSbnd} 1-\gamma= \bigl(1- \bigl(\gamma'
\bigr)^{1/\Delta} \bigr)^2>1-6^{-4}.
\end{equation}

We can now verify condition \eqref{flip2}.

%
\begin{lem}\label{lemflip2} If $\xi^\vep$ is a
voter model perturbation satisfying the hypotheses of
Theorem~\ref{thmCCTpert}, then there exists $\vep_1>0$
and $x_0\in{\mathbb{Z}^d}$ such that \eqref{flip2} holds for
$\xi^\vep$ if $\vep<\vep_1$.
\end{lem}

\begin{pf} For $\gamma'$ as above, let $\vep_1$ be as in
Lemma~\ref{vmpertspercolate}, so that for $0<\vep<\vep_1$
all the conclusions of that lemma hold, as well as the
setup \eqref{Hprop}--\eqref{percext}, with
$\rho_\infty>0$. There are two main steps in the proof.
In the first, we show that if $A\subset2{\mathbb{Z}^d}$ is large, then
for all large $n$, $\xi^\vep_{2nT}\ne\varnothing$ will imply
$V_{2n}\cap A$ is also large; see~\eqref{K0occupied}
below. To do this, we argue that there is a uniform
positive lower bound on $P_\xi(\exists z\in V_2\mbox{
with } W^{(z,2)}_m\ne\varnothing\ \forall m\ge2)$,
$\xi\notin\{\zero,\one\}$. Iteration leads to \eqref
{K0occupied}. In
the second step, we consider $A\subset{\mathbb{Z}^d}$ large, and for
$a\in A$ choose points $\ell(a)\in{\mathbb{Z}^d}$ such that $a\in
2L\ell(a)+Q(L)$. If $A$ is sufficiently large, there will
be many points $a_i\in A$ which are widely separated. By
the first step, for large $n$, there will be many points $2\ell
(a_i)\in
V_{2n}$, and for each of these there will be points
$y^0_i,y^1_i\in2L\ell(a_i)+Q(L)$ such that
$\xi^\vep_{2nT}(y^0_i)=0$ and
$\xi^\vep_{2nT}(y^1_i)=1$. Given these points, it will
follow from Lemma~\ref{vmirred} that there is a uniform
positive lower bound on the probabilities of independent
events on which $\xi^\vep_t(a_i)=1$ and
$\xi^\vep_t(a_i+x_0)=0$ for all $t\in[(2n+1)T, (2n+3)T]$.
Many of these events will occur, forcing
$A(\xi_t,x_0)$ to be large; see \eqref{emptyA}
below. Condition \eqref{flip2} now follows easily.

It is convenient to start with two estimates which
depend only on the process $\xi^\vep$ (and not on the
percolation construction). We claim that by
Lemma~\ref{vmirred} with $t=2T$,
%
%
\begin{equation}
\min_{x\in Q(L), k=1,\ldots,d}\mathop{\inf_{\xi\in\{0,1\}^{\mathbb
{Z}^d}\dvtx}}_{
\xi(x)=1,\xi(x+e_k)=0}
P_\xi\bigl(\xi^\vep_{2T}\in H \bigr)>0.
\end{equation}
To see this, note that for small $\vep$, $\xi\in
H$ depends only on the coordinates $\xi(x)$, $x\in
Q(L+1)$. This means there are disjoint sets $B_0,B_1\subset Q(L+1)$ so
that $\xi^\vep_{2T}(x)=i$ for all $x\in B_i$, $i=0,1$,
implies $\xi^\vep_{2T}\in H$. If $G(2T,y_0,y_1,B_0,B_1)$
and $\Lambda(y_0,y_1,B_0,B_1)$ are as in Lemma~\ref{vmirred}
with $(y_0,y_1)=(x+e_k,x)$, then for $x\in Q(L)$ the above
infimum is bounded below by $P(G(2T,x+e_k,x,B_0,B_1))>0$.
If $\xi\notin\{\zero,\one\}$, there must exist
$k\in\{1,\ldots,d\}$ and $x,z\in{\mathbb{Z}^d}$ with $x\in2Lz+Q(L)$ and
$\xi(x)=1$, $\xi(x+e_k)=0$. It now follows from translation
invariance that
%
%
\begin{equation}
\label{rho1} \rho_1= \inf_{\xi\notin\{\zero,\one\}}P_\xi
\bigl( \exists z\in{\mathbb{Z}^d}\mbox{ such that }
\xi^\vep_{2T}\circ\sigma_{2Lz}\in H \bigr)>0.
\end{equation}
Let $\rho_2=\rho_1\rho_\infty>0$.

Next, suppose
$y_0,y_1,y\in Q(L)$, $B_{1}=\{y\}$, $B_0=\{y+x_0\}$ and
$G(T,y_0,y_1,\break B_0, B_1)$ be as in
Lemma~\ref{vmirred}. To also
require that $\xi^\vep_u$ be constant at $y,y+x_0$ for
$u\in[T,3T]$, we let $\tilde G(T,y_0,y_1,B_0,B_1)$ be the
event
\[
G(T,y_0,y_1,B_0,B_1)\cap\bigl
\{\mbox{$T^z_m,T^{*,z}_m
\notin[T,3T]$ for $z=y,y+x_0$ and all $m\ge1$} \bigr\}.
\]
Note that each $\tilde G$ is an intersection of two independent
events each with positive probability, and so $P(\tilde G)>0$.
Making use of the notation of Lemma~\ref{vmirred}, choose $\tilde
M<\infty$ such
that
\[
\Lambda=\bigcup_{y_0,y_1,y\in Q(L)}\Lambda(y_0,y_1,B_0,B_1)
\subset[-\tilde M,\tilde M]^d,
\]
and put
%
%
\begin{equation}
\label{deltaprime} \tilde\delta=\min_{y_0,y_1,y\in Q(L)}P \bigl(\tilde
G(T,y_0,y_1,B_0,B_1) \bigr)>0.
\end{equation}
If $\xi^\vep_0(y_i)=i,i=0,1$, then
%
%
\begin{equation}
\begin{tabular}{p{260pt}@{}}
\label{Gimpl}$\tilde G(T,y_0,y_1,B_0,B_1)$
implies $\xi^\vep_t(y)=1,\xi^\vep_t(y+x_0)=0$
for all $t\in[T,3T]$.
\end{tabular}
\end{equation}

We now start the proof of
%
%
\begin{equation}
\label{0occupied2} \lim_{K\to\infty} \mathop{\sup
_{A\subset2{\mathbb{Z}^d}}}_{|A|\ge K} \lim_{n\to\infty}P_\xi
\bigl(\xi^\vep_{2nT}\ne\varnothing\mbox{ and }
V_{2n}\cap A=\varnothing\bigr) =0.
\end{equation}
Fix $\delta>0$.
By \eqref{0occupied} there exists $K_1=K_1(\delta)<\infty$
such that if $A\subset2{\mathbb{Z}^d}$ with $|A|\ge K_1$, then there exists
$n_1=n_1(A)<\infty$ such that
%
%
\begin{equation}
\label{afternA} P \bigl(W^0_{2n}\ne\varnothing,
W^0_{2n}\cap A=\varnothing\bigr) < \delta\qquad\mbox{for all
}n\ge n_1.
\end{equation}
We may increase $n_1$ if necessary so that
$P(|\cC(0,0)|<\infty, W^0_{2n_1}\ne\varnothing)<\delta$,
which implies that
%
%
\begin{equation}
\label{dieaftern12} P \bigl(W^0_{2n_1}\ne\varnothing,
W^0_{2n}=\varnothing\bigr) < \delta\qquad\mbox{for all }n\ge
n_1.
\end{equation}

For $j=1,2,\ldots,$ let $m(j)=(j-1)(2n_1+2)$,
and define a random sequence of sites $z_j$, as follows.
If $V_{m(j)+2}=\varnothing$, put $z_j=0$. If not, choose $z\in V_{m(j)+2}$
with minimal norm (with some convention for ties), and put
$z_j=z$. By the Markov property and \eqref{rho1},
%
%
\begin{equation}
\label{N1} \inf_{\xi\notin\{\zero,\one\}}P_\xi\bigl(z_1
\in V_2, \bigl|\cC(z_1,2)\bigr|=\infty\bigr) \ge\rho_2.
\end{equation}

Let
\[
N=\inf\bigl\{j\dvtx z_j\in V_{m(j)+2}\mbox{ and }
W^{z_j,m(j)+2}_{m(j+1)}\ne\varnothing\bigr\}
\]
and $\cF_n$ be the
$\sigma$-algebra generated by $\cG(\R^d\times[0,nT])$ and the
$\theta(z,k)$ for $z\in\Z^d, k<n$. It follows from our construction and
\eqref{N1} that
almost surely on the event
$\{\xi^\vep_{m(j)T}\ne\varnothing\}$,
\begin{eqnarray*}
&&P_\xi\bigl(z_j\in V_{m(j)+2} \mbox{ and }
W^{(z_j,m(j)+2)}_{m(j+1)}=\varnothing\mid\cF_{m(j)} \bigr)
\\
&&\qquad=P_{\xi^\vep_{m(j)T}} \bigl(z_1\in V_{2} \mbox{ and }
W^{(z_1,2)}_{2n_1+2}=\varnothing\bigr)
\\
&&\qquad\le P_{\xi^\vep_{m(j)T}} \bigl(z_1\in V_2, \bigl|
\cC(z_1,2)\bigr|<\infty\bigr)
\\
&&\qquad\le1-\rho_2.
\end{eqnarray*}
In the last line note that by \eqref{Hprop} if the initial state is
$\one$, the probability is zero
as~$\one$ is a trap.
Since the event on the LHS is $\cF_{m(j+1)}$-measurable, we may iterate
this inequality to obtain
%
%
\begin{equation}
\label{Ntailbound}P_\xi\bigl(\xi^{\vep}_{m(j)T}\ne
\varnothing, N>j \bigr)\le(1-\rho_2)^{j}.
\end{equation}
Taking $J_0>2$ large enough so that $(1-\rho_2)^{J_0}<\delta$,
and then $2n> m(J_0+1)$,
%
%
\begin{eqnarray}
\label{decompJ0}&& P_\xi\bigl(\xi^\vep_{2nT}\ne
\varnothing, V_{2n}\cap A=\varnothing\bigr)
\nonumber
\\[-8pt]
\\[-8pt]
\nonumber
&&\qquad\le\delta+ \sum
_{j=1}^{J_0}P_\xi\bigl(
\xi^\vep_{m(j)T}\ne\varnothing, V_{2n}\cap A=
\varnothing, N=j \bigr).
\end{eqnarray}
For $j\le J_0$, almost surely on the event
$\{\xi^\vep_{m(j)T}\ne\varnothing, N>j-1\}$,
\begin{eqnarray*}
&&P_\xi\bigl(z_j\in V_{m(j)+2},
W^{(z_j,m(j)+2)}_{m(j+1)}\ne\varnothing, V_{2n}\cap A=
\varnothing\mid\cF_{m(j)} \bigr)
\\
&&\qquad= P_{\xi^\vep_{m(j)}} \bigl(z_1\in V_{2},
W^{(z_1,2)}_{2n_1}\ne\varnothing, V_{2n-2n_1}\cap A=
\varnothing\bigr)
\\
&&\qquad\le P_{\xi^\vep_{m(j)}} \bigl(z_1\in V_{2},
W^{(z_1,2)}_{2n_1}\ne\varnothing, W^{(z_1,2)}_{2n-2n_1}
\cap A=\varnothing\bigr)
\\
&&\qquad\le\delta+P_{\xi^\vep_{m(j)}} \bigl(z_1\in V_{2},
W^{(z_1,2)}_{2n-2n_1}\ne\varnothing, W^{(z_1,2)}_{2n-2n_1}
\cap A=\varnothing\bigr)
\\
&&\qquad\le2\delta,
\end{eqnarray*}
where the last three inequalities follow from \eqref{percext},
\eqref{dieaftern12}, \eqref{afternA} and the fact that $n\ge2n_1$ by
our choice of $n$ above.
Combining this bound with \eqref{decompJ0} and then using \eqref
{Ntailbound}, we obtain
\begin{eqnarray*}
P_\xi\bigl(\xi^\vep_{2nT}\ne\varnothing,
V_{2n}\cap A=\varnothing\bigr) &\le&\delta+ 2\delta\sum
_{j=1}^{J_0}P_\xi\bigl(
\xi^\vep_{m(j)T}\ne\varnothing, N> j-1 \bigr)
\\
&\le&\delta+ 2\delta\sum_{j=1}^{J_0}(1-
\rho_2)^{j-1}
\\
&\le&\delta+2\delta/\rho_2.
\end{eqnarray*}
This
establishes \eqref{0occupied2}, which along with
the argument proving
\eqref{flip2a}, implies that
for any $K_0<\infty$,
%
%
\begin{equation}
\label{K0occupied} \lim_{K\to\infty}\mathop{\sup
_{A\subset
2{\mathbb{Z}^d}}}_{|A|\ge K} \limsup_{n\to\infty}P_\xi
\bigl(\xi^\vep_{2nT} \ne\varnothing, |V_{2n}\cap
A|\le K_0 \bigr)=0.
\end{equation}

Now fix $K_0<\infty$ so that $(1-\tilde\delta)^{K_0}<\delta$. By
\eqref
{K0occupied} there exists
$K_1<\infty$ such that for $A'\subset2\Z^d$ satisfying $|A'|\ge
K_1$, there
exists $n_1(A')$ so that
%
%
\begin{equation}
\label{end2} P_\xi\bigl(\xi^\vep_{2nT}\ne
\varnothing\mbox{ and }\bigl|V_{2n}\cap A'\bigr|\le
K_0 \bigr) <\delta\qquad\mbox{if } n\ge n_1
\bigl(A' \bigr).
\end{equation}
For $a\in{\mathbb{Z}^d}$ let $\ell(a)$ be the minimal point in some
ordering of $\Z^d$ such that $a\in2L\ell(a)+Q(L)$. For
$A\subset{\mathbb{Z}^d}$ let
$\ell(A)=\{\ell(a),a\in A\}$. With $K_0,K_1$ as above,
choose $K_2<\infty$ so that if $A\subset{\mathbb{Z}^d}$ and $|A|\ge
K_2$, then $\ell(A)$ contains $K_1$ points,
$\ell(a_1),\ldots,\ell(a_{K_1})$, such
that $|\ell(a_i)-\ell(a_{j})|2L\ge4\tilde M$ for $i\ne j$.
The regions
$2L\ell(a_1)+[-\tilde M,\tilde M]^d,\ldots,
2L\ell(a_{K_1})+[-\tilde M,\tilde M]^d$ are pairwise disjoint. Let
$A'=\{2\ell(a_1),\ldots,2\ell(a_{K_1})\}\subset2\Z^d$.

Now suppose $t\in[(2n+1)T,(2n+3)T]$ for some integer $n\ge
n_1(A')$.\break  By~\eqref{end2}, on the event
$\{|\xi^\vep_{2nT}|>0\}$, except for a set of probability at most
$\delta$, $V_{2n}$ will contain at least $K_0$ points of
$A'$. If $2\ell(a_i)$ is such a point, then by the
definitions of $V_{2n}$ and $H$, there will exist points
$y^i_0,y^i_1\in2L\ell(a_i)+Q(L)$ such that
$\xi^\vep_{2nT}(y^i_0)=0,\xi^\vep_{2nT}(y^i_1)=1$. Conditional
on this, by \eqref{deltaprime} and \eqref{Gimpl}, the
probability that
$\xi^\vep_t(a_i)=1,\xi^\vep_t(a_i+x_0)=0$ is at least
$\tilde\delta$. By independence of the Poisson point
process on disjoint space--time regions, it follows that
%
%
\begin{equation}
\label{emptyA} P \bigl(\xi^\vep_{2nT}\ne\varnothing\mbox{
and } A \bigl(x_0,\xi^\vep_t \bigr)=
\varnothing\bigr) <\delta+ (1-\tilde\delta)^{K_0},
\end{equation}
and therefore since $t>2nT$,
\[
P \bigl(\xi^\vep_{t}\ne\varnothing\mbox{ and } A
\bigl(x_0,\xi^\vep_t \bigr)=\varnothing
\bigr) <\delta+ (1-\tilde\delta)^{K_0}<2\delta,
\]
the last by our choice of $K_0$.
This proves \eqref{flip2}.

Finally, \eqref{N1} implies by
\eqref{percext}, \eqref{Hprop}, the definition of $V_n$ and the fact
that ${\one}$ is
a trap by Lemma~\ref{elemann}, that
%
%
\begin{equation}
\label{survival} \inf_{\xi\neq\zero}P_\xi\bigl(
\xi^\vep_t\ne\varnothing\ \forall t\ge0 \bigr)\ge
\rho_2.
\end{equation}
This will be used below.
\end{pf}

\begin{pf*}{Proof of Theorem~\ref{thmCCTpert}}
We verify the assumptions of
Proposition~\ref{propCCT}. It follows from Lemma~\ref{elemann} and
\eqref{traps} that $c_\vep(x,\xi)$ is symmetric and
$\zeta^\vep_t$, the annihilating dual of $\xi^\vep_t$, is
parity preserving. By Corollary~\ref{vpirred} (which
applies by Remark~\ref{nvm}) there exists $\vep_3>0$ such
that if $0<\vep<\vep_3$, then $\zeta^\vep_t$ is
irreducible. By Lemmas~\ref{lemodd2} and \ref{lemflip2}
(and the proof of the latter), there exists
$0<\vep_4<\vep_3$ such that if $0<\vep<\vep_4$, then
\eqref{flip2}, \eqref{oddgoal} and \eqref{survival} hold
for $\xi^\vep_t$.

Assume now that $0<\vep<\vep_4$. It remains to
check that the dual growth condition~\eqref{eqH} (the
conclusion of Lemma~\ref{lemH}) holds, and to do this
it suffices by Remark~\ref{remHaltcond} to show that
\eqref{eqH3} for $\xi^\vep$ holds. By \eqref{flip2} and
\eqref{survival}, there is a $\delta_1>0$, $t_0<\infty$ and
$A\in Y$ so that for all $t\ge t_0$ (with
$\xi^\vep_0=1_{\{0\}}$),
\[
P \bigl(\xi^\vep_t(a)=1\mbox{ for some }a\in A \bigr)
\ge P \bigl(A \bigl(x_0,\xi_t^\vep\bigr)
\neq\varnothing\bigr)\ge\delta_1.
\]
Next
apply \eqref{02z}, translation invariance and the Markov
property to conclude that for $t$ as above,
\begin{eqnarray*}
P \bigl(\xi^\vep_{t+1}(0)=1 \bigr)&\ge& E \bigl(1 \bigl(
\xi^\vep_t(a)=1\mbox{ for some }a\in A
\bigr)P_{\xi^\vep_t} \bigl(\xi_1^\vep(0)=1 \bigr) \bigr)
\\
&\ge&\delta_1\min_{a\in A}\inf_{\xi^\vep_0\dvtx\xi^\vep
_0(0)=1}P_{\xi_0^\vep
}
\bigl(\xi_1^\vep(-a)=1 \bigr)\ge\delta_2>0.
\end{eqnarray*}
This proves \eqref{eqH3}, and all the assumptions of
Proposition~\ref{propCCT} have now been verified for
$\xi^\vep_t$ if $0<\vep<\vep_4$, and thus the weak limit
\eqref{eqgoal} also holds. Finally, by~\eqref{betanonzero}
this result implies the full complete convergence\vspace*{1pt} theorem
with coexistence if $|\hat\xi_0^\vep|=\infty$. If $|\hat
\xi_0^\vep|<\infty$, then $|\xi^\vep_0|=\infty$, and the
result now follows by the symmetry of $\xi^\vep$; recall
Lemma~\ref{elemann}.
\end{pf*}

\section{\texorpdfstring{Proof of Theorem~\protect\ref{thmLVCCT}}
{Proof of Theorem 1.1}}\label{secalmostlast}
Let us check that $\mathrm{LV}(\alpha)$, $\alpha\in(0,1)$, is cancellative.
(This was
done in \cite{NP} for the case $p(x)=1_{\cN}(x)/|\cN|$ for
$\cN$ satisfying \eqref{N}.) For the more general setting here, we
assume $p(x)$ satisfies
\eqref{passump}, and allow any $d\ge1$. We first observe
that if $c(x,\xi)$ has the form given in \eqref{q00}, then
it follows from \eqref{q01} that
\[
c(0,\xi) = k_0\sum_{A\in Y}
q_0(A) \frac12 \bigl[1 - \bigl(2\xi(0)-1 \bigr)H(\xi,A) \bigr].
\]
From this it is clear that the sum of two positive multiples of
cancellative rate functions is cancellative.
It follows from a bit of arithmetic that if \eqref{diag} holds, then
$\mathrm{LV}(\alpha)$ with $\vep^2=1-\alpha>0$ has
flip rates
\[
c_{\mathrm{LV}}(x,\xi)=\alpha c_{\vm}(x,\xi)+\vep^2f_0f_1(x,
\xi).
\]
We have already noted that $c_{\vm}$ is cancellative, and so by the
above we
need only check that $c^*(x,\xi)=f_0(x,\xi)f_1(x,\xi)$ is
cancellative.

To do this we let $p^{(2)}(0)=\sum_{x\in{\mathbb{Z}^d}}(p(x))^2$,
$k_0 = (1-p^{(2)}(0))/2$, $q_0(A)=0$ if $|A|\ne3$
and
\[
q_0 \bigl(\{0,x,y\} \bigr)= k_0^{-1}p(x)p(y)\qquad
\mbox{if $0,x,y$ are distinct}.
\]
Note that $\sum_{A\in Y}q_0(A)=1$
because [recall that $p(0)=0$]
\[
\sum_{\{x,y\}}q_0 \bigl(\{0,x,y\} \bigr) =
\frac{1}{2k_0}\sum_{x\ne y}p(x)p(y) =
\frac{1}{2k_0} \bigl(1-p^{(2)}(0) \bigr) = 1.
\]
Also, for $0,x,y$ distinct,
\[
\tfrac12 \bigl[1- \bigl(2\xi(0)-1 \bigr)H \bigl(\xi,\{0,x,y\} \bigr)
\bigr] = 1
\bigl\{ \xi(x)\ne\xi(y) \bigr\}.
\]
With these facts it is easy to see that
\begin{eqnarray*}
&&k_0\sum_{A\in Y} q_0(A)
\frac12 \bigl[1 - \bigl(2\xi(0)-1 \bigr)H(\xi,A) \bigr]
\\
&&\qquad= \sum_{x,y} p(x)p(y)1 \bigl\{\xi(x)\ne\xi(y) \bigr
\} = f_0(0,\xi)f_1(0,\xi),
\end{eqnarray*}
proving $c^*(x,\xi)=f_0(x,\xi)f_1(x,\xi)$ is cancellative and hence so
is $\mathrm{LV}(\alpha)$.

Although we won't need it, we calculate the parameters of the branching
annihilating dual. Adding in the voter model, we see that they are
\begin{eqnarray*}
k_0&=& \alpha+ (1-\alpha)\frac{1-p^2(0)}{2},\qquad q_0 \bigl(
\{y\} \bigr) = \frac{\alpha}{k_0} p(y),\\
 q_0 \bigl(\{0,y,z\} \bigr)& =&
\frac{1-\alpha}{k_0} p(y)p(z),
\end{eqnarray*}
and $q_0(A)=0$ otherwise. One can see from this that
$\zeta_t$, the dual of $\mathrm{LV}(\alpha)$, describes a system of
particles evolving according to the following rules:
(i) a particle at $x$ jumps to $y$ at rate $\alpha p(y-x)$;
(ii) a particle at $x$ creates two particles and sends them to
$y,z$ at rate $(1-\alpha) p(y-x)p(z-x)$; (iii)
if a particle attempts to land on another particle, then
the two particles annihilate each other.

Assume $d\ge3$. The function $f(u)$ as shown in Section~1.3 of
\cite{CDP11} is a cubic, and under the assumption
\eqref{diag} reduces to $f(u)=2p_3(1-\alpha)u(1-u)(1-2u)$, where
$p_3$ is a certain (positive) coalescing random walk
probability. Thus $f'(0)>0$, so the complete convergence
theorem with coexistence for $\mathrm{LV}(\alpha)$ for $\alpha$
sufficiently close to one follows from
Theorem~\ref{thmCCTpert}.

Now suppose $d=2$. It suffices to prove an analogue of
Lemma~\ref{vmpertspercolate} as the above results will then
allow us to apply the proof of Theorem~\ref{thmCCTpert} in
the previous section to give the result. As
the results of \cite{CDP11} do not apply, we will use results
from \cite{CMP} instead and proceed as in Section~4 of
\cite{CP07}. Instead of \eqref{exptail}, we only require
(as was the case in \cite{CMP})
%
%
\begin{equation}
\label{3mom} \sum_{x\in\Z^2}|x|^3p(x)<
\infty.
\end{equation}
We will need some notation from \cite{CMP}. For $N>1$, let $\xi^{(N)}$
be the $\mathrm{LV}(\alpha_N)$ process where
\[
\alpha_N=1-\frac{(\log N)^3}{N},
\]
and consider the rescaled process,
$\xi^N_t(x)=\xi^{(N)}_{Nt}(x\sqrt N)$, for $x\in S_N=\Z^2/\sqrt N$. The
associated process taking values in $M_F(\R^2)$ (the space of
finite\vadjust{\goodbreak}
measures on the plane with the weak topology) is
%
%
\begin{equation}
\label{MVP} X_t^N=\frac{\log N}{N}\sum
_{x\in S_N}\xi_t^N(x)
\delta_x.
\end{equation}
For parameters $K,L'\in\NN$, $K>2$ and $L'>3$, which will be chosen
below, we let $\underline\xi^N_t(x)\le\xi^N_t(x)$, $x\in S_N$ be a
coupled particle system where particles are ``killed'' when they exit
$(-KL',KL')^2$, as described in Proposition~2.1 of \cite{CP07}. (Here a
particle corresponds to a $1$.) In particular $\underline\xi^N_t(x)=0$
for all $|x|\ge KL'$. $\underline X^N_t$ is defined as in \eqref{MVP}
with $\underline\xi^N$ in place of $\xi^N$.

We will need to keep track of some of the dependencies in the constant
$C_{8.1}$ in Lemma~8.1 of \cite{CMP}. As in that result, $B^N$ is a
rate $N\alpha_N=N-(\log N)^3$ random walk on $S_N$ with step
distribution $p_N(x)=p(x\sqrt N)$, $x\in S_N$, starting at the origin.

%
\begin{lem}\label{lem81mod} There are positive constants
$c_0$ and $\delta_0$ and a nondecreasing function
$C_0(\cdot)$, so that if $t>0$, $K,L'\in\NN$, $K>2$ and
$L'>3$, and $X_0^N=\underline X_0^N$ is supported on $[-L',L']^2$, then
%
%
\begin{eqnarray}
\label{killbnd} &&E \bigl(X_t^N(1)-\underline
X_t^N(1) \bigr)
\nonumber
\\
&&\qquad\le X_0^N(1)
\Bigl[c_0e^{c_0 t}P\Bigl(\sup_{s\le t}
\bigl|B^{N}_s\bigr|>(K-1)L'-3 \Bigr)\\
&&\hspace*{52pt}\qquad\quad{}+C_0(t) \bigl(1\vee X_0^N(1)
\bigr) (\log N)^{-\delta_0} \Bigr].\nonumber
\end{eqnarray}
\end{lem}
\begin{pf} This is a simple matter of keeping track of the
$t$-dependency in some of the constants arising in the proof of
Lemma~8.1 in \cite{CMP}.
\end{pf}

Recall from Theorem~1.5 of \cite{CMP} that if $X_0^N\to X_0$ in
$M_F(\R
^2)$, then $\{X^N\}$ converges weakly in $D(\R_+,M_F(\R^2))$ to a
two-dimensional super-Brownian motion, $X$, with branching rate $4\pi
\sigma^2$, diffusion coefficient $\sigma^2$ and drift $\eta>0$ [write
$X$ is $\mathrm{SBM}(4\pi\sigma^2,\sigma^2,\eta)$], where $\eta$ is the constant
$K$ in (6) of \cite{CMP} (not to be confused with our parameter $K$).
See (MP) in Section~1 of \cite{CMP} for a precise definition of SBM.
The important point for us is that the positivity of $\eta$ will mean
that the supercritical $X$ will survive with positive probability, and
on this set will grow exponentially fast.

We next prove a version of Proposition~4.2 of \cite{CP07} which when
symmetrized is essentially a scaled version of the required Lemma~\ref
{vmpertspercolate}. To be able to choose $\gamma'$ as in \eqref{LSSbnd},
so that we may apply Lemma~5.1 of \cite{CDP11}, we will have to be more
careful with the selection of constants in the proof of Proposition~4.2
in the above reference. We start by choosing $c_1>0$ so that
%
%
\begin{equation}
\label{c1choice} \bigl(1-e^{-c_1} \bigr)^2>1-6^{-4},
\end{equation}
and then setting
\[
\gamma'_K=e^{-c_1(2K+1)^3}.\vadjust{\goodbreak}
\]

%
\begin{lem}\label{lemprop42} There are $T'>1$, $L',K,J'\in\NN$ with
$K>2$, $L'>3$, and $\vep_1\in(0,\frac{1}{2})$ such that if
$0<1-\alpha
<\vep_1$, $N>1$ is chosen so that $\alpha=1-\frac{(\log N)^3}{N}$, and
$I_{\pm e_i}=\pm2L'e_i+[-L',L']^2$, then
\[
\underline X_0^N \bigl( \bigl[-L',L'
\bigr]^2 \bigr)\ge J'\mbox{ implies }P \bigl(\underline
X_{T'}^N(I_{\pm e_i})\ge J' \mbox{
for }i=1,2 \bigr)\ge1-\gamma'_K.
\]
\end{lem}
\begin{pf} By the monotonicity of $\underline X^N$ in its
initial condition (Proposition~2.1(b) of \cite{CP07} and
the monotonicity of $\mathrm{LV}(\alpha)$ discussed, e.g.,
in Section 1 of \cite{CP07}), we may assume that
$\underline X_0^N(\R^2\setminus[-L',L']^2)=0$ and
$\underline X_0^N([-L',L']^2)\in[J',2J']$, where $L'$ and $J'$ are
chosen below.

We will choose a number of constants which depend on an integer
$K>2$ and will then choose $K$ large enough near the end
of the proof. Assume $B=(B^1,B^2)$ is a $2$-dimensional
Brownian motion with diffusion parameter $\sigma^2$,
starting at $x$ under $P_x$ and fix $p>\frac{1}{2}$. Set
%
%
\begin{equation}
\label{Tdef} T'=c_2K^{2p},
\end{equation}
where a short calculation shows that if $c_2$ is chosen large enough,
depending on $\sigma^2$ and $\eta$, then for any $K>2$,
%
%
\begin{equation}
\label{Tdef2} e^{\eta T'/2}\inf_{|x|\le K^p} P_x
\bigl(B_1\in\bigl[K^p,3K^p
\bigr]^2 \bigr)\ge5.
\end{equation}
Now put $L'=K^p\sqrt{T'}$, increasing $c_2$ slightly so that
$L'\in\NN$. If $I=[-L',L']^2$ and $X$ is the limiting super-Brownian motion
described above, then as in Lemma~12.1(b) of \cite{DP99},
there is a $c_3(K)$ so that
%
%
\begin{equation}
\label{SBMbnd}\qquad \forall J'\in\NN\mbox{ and }i\le2,\mbox{ if
}X_0(I)\ge J',\mbox{ then }P \bigl(X_{T'}(I_{\pm e_i})<4J'
\bigr)\le c_3/J'.
\end{equation}
Next choose $J'=J'(K)\in\NN$ so that
\[
\frac{c_3}{J'}\le\frac{\gamma'_K}{100}.
\]
As in Lemma~4.4 of \cite{CP07}, the weak convergence of $X^N$ to $X$
and \eqref{SBMbnd} show that for $N\ge N_1(K)$,
%
%
\begin{equation}
\label{massbnd} \forall i\le2\mbox{ if }X_0^N(I)\ge
J', \mbox{ then }P \bigl(X_{T'}^N(I_{\pm
e_i})<4J'
\bigr)\le\frac{\gamma'_K}{50}.
\end{equation}

Next use Lemma~\ref{lem81mod}, the fact that
$X^N_{T'}-\underline X^N_{T'}$ is a nonnegative measure
and Donsker's theorem to see that there is a $c_4>0$ and an
$\vep_N=\vep_N(K)\to0$ as $N\to\infty$, so that for any
$i\le2$,
\begin{eqnarray*}
&&P \bigl(X_{T'}^N(I_{\pm e_i})-\underline
X_{T'}^N(I_{\pm e_i})\ge2J' \bigr)
\\
&&\qquad\le\frac{X_0^N(1)}{2J'} \Bigl[c_0e^{c_0T'}
\Bigl(P_0 \Bigl(\sup_{s\le
T'}|B_s|>(K-1)L'-3
\Bigr)+\vep_N \Bigr)
\\
&&\hspace*{96pt}\qquad\quad{} +C_0 \bigl(T' \bigr) \bigl(1\vee
X_0^N(1) \bigr) (\log N)^{-\delta_0} \Bigr]
\\
&&\qquad\le\bigl[c'_0e^{c_0T'} \bigl(\exp
\bigl(-c_4K^{2+2p} \bigr)+\vep_N
\bigr)+C_0 \bigl(T' \bigr)2J'(\log
N)^{-\delta_0} \bigr],
\end{eqnarray*}
where the fact that $X_0^N(1)\le2J'$ and the definition of $L'$ are
used in the last line. It follows that for $K\ge K_0$ and $N\ge
N_2(K)$, the above is bounded by
%
%
\begin{equation}
\label{survbnd} 2c'_0e^{c_0T'}\exp
\bigl(-c_4K^{2+2p} \bigr)\le2c'_0e^{c_5K^{2p}-c_4K^{2+2p}}
\le\frac{\gamma'_K}{50}.
\end{equation}
The fact that $p>\frac{1}{2}$ is used in the last inequality.
We finally choose $K\in\NN^{>2}$, $K\ge K_0$. Therefore the bounds in
\eqref{massbnd} and \eqref{survbnd} show that for $N\ge N_1(K)\vee
N_2(K)$ and $i\le2$,
\begin{eqnarray*}
P \bigl(\underline X_{T'}^N(I_{\pm e_i})<2J'
\bigr)&\le& P \bigl(X^N_{T'}(I_{\pm e_i})
\le4J' \bigr)+P \bigl(X_{T'}^N(I_{\pm e_i})-
\underline X_{T'}^N(I_{\pm e_i})
\ge2J' \bigr)
\\
&\le&\frac{\gamma'_K}{25}.
\end{eqnarray*}
Sum over the $4$ choices of $\pm e_i$ to prove the required result
because the condition on $N$ is implied by taking $1-\alpha=(\log
N)^3/N$ small enough.
\end{pf}
\begin{pf*}{Completion of proof of Theorem~\ref{thmLVCCT}} By
symmetry we have an analogue of the above lemma with $0$'s in place of
$1$'s. Let $\alpha$ and $N$ be as in Lemma~\ref{lemprop42}. Now undo
the scaling and set $L=\sqrt N L'$, $J=\frac{N}{\log N}J'$ and $T=T'
N$. Slightly abusing our earlier notation we let $\underline{\xi
}_t\le
\xi^{(N)}_t$ be the unscaled coupled particle system
where particles are killed upon exiting $(-KL,KL)^2$ and let $\tilde
I_{\pm e_i}=\pm e_iL+[-L,L]^2$. We define
\[
G_\xi= \bigl\{\underline{\xi}_T(\tilde
I_{\pm e_i})\ge J, \hat{\underline\xi}_T(\tilde
I_{\pm e_i})\ge J\mbox{ for }i=1,2 \bigr\},
\]
where ${\underline\xi}_0=\xi$.
Lemma~\ref{lemprop42} gives the conclusion of Lemma~\ref
{vmpertspercolate} with $\vep=1-\alpha$, $\gamma'=\gamma'_K$ and
now with
\[
H= \bigl\{\xi\in\{0,1\}^{{\mathbb{Z}^d}}\dvtx\xi\bigl([-L,L]^d \bigr)
\ge J, \hat\xi\bigl([-L,L]^d \bigr)\ge J \bigr\}.
\]
Note that by \eqref{c1choice} and the definition of $\gamma'_K$, we
have $1-\gamma>1-6^{-4}$ where
$\gamma$ is as in \eqref{LSSbnd}. The definition of
$\underline\xi$ gives the required measurability of
$G_\xi$. Note that $H$ depends only on $\{\xi(x)\dvtx x\in[-L,L]^d\}
$, and
$\xi\in H$ implies $\xi(x)=1$ and $\xi(x')=0$ for some $x,x'\in
[-L,L]^d$. These are the only properties of $H$ used in the previous
proof. Finally it is easy to adjust the parameters so that $L\in\NN$ as
in Lemma~\ref{vmpertspercolate}. One way to do this is to modify
\eqref
{SBMbnd} so the conclusion of Lemma 6.2 becomes
\[
\underline X_0^N \bigl(I' \bigr)\ge
J'\mbox{ implies }P \bigl(\underline X_T^N
\bigl(I'_{\pm
e_i} \bigr)\ge J' \mbox{ for
}i=1,2 \bigr)\ge1-\gamma'_K,
\]
where $I'=[-L'-1,L'+1]^2$ and $I'_{\pm e_i}=\pm2L'e_i+[-L'+1,L'-1]^2$.
Then for $N$ large (in addition to the constraints above, $N\ge9$ will
do) one can easily check that the above argument is valid with
$L=\lfloor\sqrt{N}L'\rfloor\in\NN$.
Therefore,
with the conclusion of this version of Lemma~\ref{vmpertspercolate} in
hand, the result for $d=2$ now follows as in the proof of Theorem~\ref
{thmCCTpert}.
\end{pf*}

%
\begin{rem}\label{noexptail} The above argument works equally well
for $\mathrm{LV}(\alpha)$ for $d\ge3$ even without assuming \eqref{3mom}. Only
a few constants need to be altered, for example, $p=(d-1)/2$ and
$\gamma
'_K=e^{c_1(2K+1)^{d+1}}$. More generally the argument is easily
adjusted to give the result
for the general voter model perturbations in Theorem~\ref{thmCCTpert}
(for $d\ge3$) without assuming \eqref{exptail}, provided the particle
systems are also attractive. This last condition is needed to use the
results in \cite{CP07}.
\end{rem}

\section{\texorpdfstring{Proofs of Theorems~\protect\ref{thmtvm} and \protect\ref{thmgeom}}
{Proofs of Theorems 1.3 and 1.4}}\label{seclast}
\mbox{}
\begin{pf*}{Proof of Theorem~\ref{thmtvm}}
Let $\xi_t$ be the affine voter model with parameter
$\alpha\in(0,1)$, and $d\ge3$.
If $\vep^2=1-\alpha$, then the rate function of $\xi$ is of the form in~\eqref{vmpert} and
\eqref{vmpert2} where
%
%
\begin{equation}
\label{avh} h_i(x,\xi)=-f_i(x,\xi)+1 \bigl(\xi(y)=i
\mbox{ for some }y\in\cN\bigr).
\end{equation}
Taking $Z_1,\ldots,Z_{N_0}$ to be the distinct points in $\cN$ we see
that $\mathrm{AV}(\alpha)$ is a voter model perturbation. The fact that
$c_{\tvm
}(x,\xi)$ is
cancellative was established in Section~2 of \cite{CD91},
and so, as for $\mathrm{LV}(\alpha)$, we may conclude
that $\mathrm{AV}(\alpha)$ is a cancellative process. It is easy to check that
$c_{\AV}(x,\xi)$ is not a pure voter
model rate function, so the only remaining condition of Theorem~\ref
{thmCCTpert} to
check is $f'(0)>0$.

To compute $f(u)$, let $\{B^x_u,u\ge0, x\in{\mathbb{Z}^d}\}$ be a
system of coalescing random walks with step distribution
$p(x)$, and put $A^F_t=\{B^x_t,x\in F\}$, $F\in Y$. The
slight abuse of notation
$|A^F_\infty|=\lim_{t\to\infty}|A^F_t|$ is
convenient. If $\xi_0(x)$ are i.i.d.
Bernoulli with $E(\xi_0(x))=u$, and $F_0,F_1\in Y$ are
disjoint, then (see (1.26) in~\cite{CDP11})
%
%
\begin{eqnarray}
\label{coaldual2} &&\bigl\langle\xi(y)= 0\ \forall y\in F_0, \xi(x)=
1\ \forall x\in F_1 \bigr\rangle_u
\nonumber
\\[-8pt]
\\[-8pt]
\nonumber
&&\qquad= \sum_{i,j}(1-u)^iu^j P
\bigl(\bigl|A^{F_0}_\infty\bigr|=i, \bigl|A^{F_1}_\infty\bigr|=j,
\bigl|A^{F_0\cup F_1}_\infty\bigr|=i+j \bigr).
\end{eqnarray}

From \eqref{f} and \eqref{avh} we have $f(u)=G_0(u)-G_1(u)$, where
\begin{eqnarray*}
G_0(u) &=& \bigl\langle1 \bigl\{\xi(0)=0 \bigr\}
\bigl(-f_1(0, \xi)+ 1 \bigl\{\xi(y)\ne0\mbox{ for some }y\in\cN\bigr
\} \bigr) \bigr\rangle_u,
\\
G_1(u) &=& \bigl\langle1 \bigl\{\xi(0)=1 \bigr\}
\bigl(-f_0(0, \xi)+ 1 \bigl\{\xi(y)\ne1\mbox{ for some }y\in\cN\bigr
\} \bigr) \bigr\rangle_u.
\end{eqnarray*}
If $c_0=\sum_ep(e)P(|A_\infty^{\{0,e\}}|=2)$, then the assumption that
$0\notin\cN$ and \eqref{coaldual2}
imply
\begin{eqnarray*}
G_0(u)&=&-c_0u(1-u)+ \bigl\langle1 \bigl\{\xi(0)=0
\bigr\} \bigr\rangle_u - \bigl\langle1 \bigl\{\xi(0)=\xi(y)=0\mbox{
for all }y\in\cN\bigr\} \bigr\rangle_u
\\
&=&-c_0u(1-u)+1-u - \sum_{j=1}^{|\cN|+1}
(1-u)^jP \bigl(\bigl|A_\infty^{\cN\cup\{0\}}\bigr|=j \bigr).
\end{eqnarray*}
Similarly,
\[
G_1(u) = -c_0u(1-u)+ u - \sum
_{j=1}^{|\cN|+1} u^jP \bigl(\bigl|A_\infty^{\cN\cup\{0\}}\bigr|=j
\bigr).
\]
Therefore if $A=|A_\infty^{\cN\cup\{0\}}|$, we obtain
\begin{eqnarray*}
f'(0) &=& G_0'(0)-G_1'(0)=-1+
\sum_{j=1}^{|\cN|+1}jP(A=j)-1+P(A=1)
\\
&=&E \bigl(A-1-1(A>1) \bigr).
\end{eqnarray*}
Note that since $A$ is $\NN$-valued, we have $A-1-1(A>1)\ge
0$ with equality holding if and only if $A\in\{1,2\}$. Hence to show
$f'(0)>0$ it suffices to establish that $P(A>2)>0$. But
since $|\cN\cup\{0\}|\ge3$ by the symmetry assumption on
$\cN$, the required inequality is easy to see by the
transience of the random walks $B^x_u$. The complete
convergence theorem with coexistence holds if $\vep>0$ is
small enough, depending on $\cN$, by
Theorem~\ref{thmCCTpert}.
\end{pf*}

\begin{pf*}{Proof of Theorem~\ref{thmgeom}}
Let $\eta^\theta_t$ be the geometric voter model with rate
function given in~\eqref{geom1}. Then $\eta^\theta_t$ is
cancellative for all $\theta\in[0,1]$ (see Section~2 of~\cite{CD91}),
and it is clear that $\eta^\theta_t$ is not a pure voter
model for $\theta<1$. [The
latter follows from the fact that $q_0(A)>0$ for any odd
subset of $\cN\cup\{0\}$.] The next step is to check that
$\eta^\theta_t$ is a voter model perturbation. Clearly $\mathbf{0}$ is a trap.
If we set $\vep^2=1-\theta$ and
$a_j=c(0,\xi)$ for $\xi(0)=0$ and
$\sum_{x\in\cN}\xi(x)=j$, then
\begin{eqnarray*}
a_j&=& \Biggl[\sum_{k=1}^j
\pmatrix{j\cr k} \bigl(-\vep^2 \bigr)^k \Biggr] \Bigg/ \Biggl[
\sum_{k=1}^{|\cN|} \pmatrix{|\cN|\cr k} \bigl(-
\vep^2 \bigr)^k \Biggr]
\\
&=&\frac{j\vep^2-{j\choose2}\vep^4+O(\vep^6)}{|\cN|\vep^2-{|\cN
|\choose2}\vep^4+O(\vep^6)},
\end{eqnarray*}
where ${j\choose2}=0$ if $j=1$.
A straightforward
calculation [we emphasize
that $c_{\vm}$ and $f_0,f_1$ are defined using
$p(x)=1_{{\cN}}(x)/|\cN|$ which satisfies
\eqref{passump} and \eqref{exptail}] now shows that
%
%
\begin{equation}
\label{geompert}\qquad c_{\gv}(x,\xi) = c_{\vm} (x,\xi) +
\vep^2\frac{|\cN|}{2} f_0(x,\xi)f_1(x,\xi)
+O \bigl(\vep^4 \bigr)\qquad\mbox{as }\vep\to0,
\end{equation}
where the $O(\vep^4)$ term is uniform in $\xi$ and may be written as a
function of $f_1(0,\xi)$.
It follows from Proposition 1.1 of \cite{CDP11} and symmetry that $\xi
^\vep$ is a voter model perturbation.

To apply Theorem~\ref{thmCCTpert} it only remains to check that
$f'(0)>0$, where
\begin{eqnarray*}
f(u) &=& \bigl\langle \bigl(1-2\xi(0) \bigr)f_1(0,
\xi)f_0(0,\xi) \bigr\rangle_u
\\
&=&\sum_{x,y}p(x)p(y) \bigl\langle \bigl(1-2\xi(0)
\bigr)\xi(x) \bigl(1-\xi(y) \bigr) \bigr\rangle_u.
\end{eqnarray*}
Using \eqref{coaldual2} it is easy to see that for $x,y,0$ distinct,
\begin{eqnarray*}
\bigl\langle\xi(x) \bigl(1-\xi(y) \bigr) \bigr\rangle_u &=&u(1-u)P
\bigl(\bigl|A^{\{x,y\}}_\infty\bigr|=2 \bigr),
\\
\bigl\langle\xi(0)\xi(x) \bigl(1-\xi(y) \bigr) \bigr\rangle_u
&=&u(1-u)P \bigl(\bigl|A^{\{0,x\}}_\infty\bigr|=1,\bigl|A^{\{0,x,y\}}_\infty\bigr|=2
\bigr)
\\
&&{} +u^2(1-u)P \bigl(\bigl|A^{\{0,x,y\}}_\infty\bigr|=3 \bigr).
\end{eqnarray*}
If we plug the decomposition ($x,y,0$ still distinct)
\begin{eqnarray*}
P \bigl(\bigl|A^{\{x,y\}}_\infty\bigr|&=&2 \bigr)= P \bigl(\bigl|A^{\{0,x\}}_\infty
\bigr|=1,\bigl|A^{\{x,y\}}_\infty\bigr|=2 \bigr)
\\
&&{}+ P \bigl(\bigl|A^{\{0,y\}}_\infty\bigr|=1,\bigl|A^{\{x,y\}}_\infty\bigr|=2
\bigr) + P \bigl(\bigl|A^{\{0,x,y\}}_\infty\bigr|=3 \bigr)
\end{eqnarray*}
into the above we find that
\[
f(u) = u(1-u) (1-2u)\sum_{x,y}p(x)p(y) P
\bigl(\bigl|A^{\{0,x,y\}}_\infty\bigr|=3 \bigr),
\]
and thus $f'(0)=\sum_{x,y}p(x)p(y)P(|A^{\{0,x,y\}}_\infty|=3)>0$
as required.
\end{pf*}

\section*{Acknowledgment}
We thank an anonymous referee for a
careful reading of the paper and making a number of helpful suggestions
which have improved the presentation.

%


\printaddresses

\end{document}